\documentclass[a4paper,twoside]{article}      
\usepackage{amsmath,amssymb,amsfonts,amsthm,amscd}
\usepackage{graphicx}
\usepackage{a4wide}
\usepackage{enumerate}
\usepackage{color}
\usepackage{pdfpages}
\usepackage[authoryear]{natbib}
\usepackage{authblk}

\newtheorem{theorem}{Theorem}[section]

\theoremstyle{definition}
\newtheorem{remark}[theorem]{Remark}
\newtheorem{algorithm}[theorem]{Algorithm}

\newtheorem{example}{Example}[section]

\usepackage{hyperref}
\hypersetup{
    bookmarks=true,         
    unicode=false,          
    pdftoolbar=true,        
    pdfmenubar=true,        
    pdffitwindow=false,     
    pdfstartview={FitH},    
    pdftitle={My title},    
    pdfauthor={Author},     
    pdfsubject={Subject},   
    pdfcreator={Creator},   
    pdfproducer={Producer}, 
    pdfkeywords={keywords}, 
    pdfnewwindow=true,      
    colorlinks=true,       
    linkcolor=blue,          
    citecolor=blue,        
    filecolor=magenta,      
    urlcolor=cyan           
}

\author[1]{Manh Hong Duong}
\affil[1]{School of Mathematics, University of Birmingham, Birmingham B15 2TT, UK. Email: h.duong@bham.ac.uk}
\author[2]{The Anh Han}
\affil[2]{School of Computing, Engineering and Digital Technologies, Teesside University, TS1 3BX, UK. Email: T.Han@tees.ac.uk}
\title{Cost efficiency of institutional incentives in finite populations}
\date{}
\begin{document}
\maketitle

\section*{Abstract}
Institutions can  provide incentives to  increase  cooperation behaviour in a population where this behaviour is infrequent.
This process is costly, and it is thus important to optimize the overall spending. 
This problem can be mathematically formulated as a multi-objective optimization problem where one wishes to minimize the cost of providing incentives while ensuring a desired level of cooperation within the population. In this paper, we provide a rigorous analysis for this problem.
 We study cooperation dilemmas in both the pairwise (the Donation game) and multi-player  (the Public Goods game) settings. 
 We prove the regularity of the (total incentive) cost function, characterize its asymptotic limits (infinite population, weak selection and large selection) and show exactly when reward or punishment is more efficient.  
 We prove that the cost function exhibits a phase transition phenomena when the intensity of selection varies.  
 We calculate the critical threshold in regards to the phase transition and study the optimization problem when the intensity of selection is under and above the critical value. It allows us to provide an exact calculation for the optimal cost of incentive, for a given intensity of selection. 
Finally, we provide numerical simulations to demonstrate the analytical results.
 Overall, our analysis provides for the first time a  selection-dependent  calculation  of the optimal cost of institutional incentives (for both reward and punishment)  that guarantees a minimum amount of cooperation. It is of crucial importance for  real-world applications of institutional incentives since intensity of selection is specific to a given population and the underlying game payoff structure.   \\

\noindent \textbf{Key words}: institution, incentives, cost optimization, evolutionary game theory, evolution of cooperation.


\newpage
\section{Introduction}
The problem of promoting the evolution of cooperative behaviour within populations of self-regarding individuals has been intensively investigated across diverse fields of behavioural, social and computational sciences \citep{perc2017statistical,key:novak2006,HanBook2013,West2007,key:Sigmund_selfishnes}. 
Various mechanisms responsible for promoting  the emergence and stability of cooperative behaviours among such individuals have been proposed. They include   kin and group selection \citep{Hamilton1964,traulsen:2006:pnas},  direct and indirect reciprocities  ~\citep{Ohtsuki2006435,krellner2020putting,nowak:2005:nature,key:HanetalAlife,okada2020review},  spatial networks  \citep{santos2006pnas,perc2013evolutionary,antonioni2017coevolution,pena2016evolutionary}, reward and punishment  ~\citep{fehr2000cooperation,boydPNAS2003,Sigmund2001PNAS,Herrmann2008,key:Hauert2007,boyd2010coordinated},  and pre-commitments \citep{nesse2001evolution,han2013good,martinez2017agreement,HanJaamas2016,sasaki2015commitment}. 
Institutional incentives, namely,  rewards for cooperation and punishment of wrongdoing, are among the most important ones \citep{wang2019exploring,Sigmund2001PNAS,han2018cost,sigmundinstitutions,vasconcelos2013bottom,chen2015first,wu2014role,garcia2019evolution,gois2019reward,powers2018modelling}. Differently from other mechanisms, in order to carry out institutional  incentives,  it is  assumed that there exists an external decision maker (e.g. institutions such as the United Nations and the European Union) that has a budget to interfere in the population to achieve a desirable outcome. Institutional enforcement mechanisms are crucial for enabling large-scale cooperation. Most modern societies implemented different forms of institutions for governing and promoting collective behaviors, including cooperation, coordination, and technology innovation \citep{ostrom1990governing,bowles2009microeconomics,bowles2002social,bardhan2005institutions,han2020Incentive,scotchmer2004innovation}.

Providing incentives is costly and it is important to minimize the cost while ensuring a sufficient level of cooperation \citep{ostrom1990governing,chen2015first,han2018cost}. 
Despite of its paramount importance, so far there have been only few works exploring this question. In particular, \cite{wang2019exploring} use optimal control theory to provide an analytical solution for  cost optimization of institutional incentives assuming deterministic evolution and infinite population sizes (modeled using replicator dynamics). This work therefore does not take into account various stochastic effects of evolutionary dynamics such as mutation and non-deterministic behavioral update \citep{traulsen2006,key:Sigmund_selfishnes,key:Hofbauer1998}. Defective behaviour can reoccur over time via mutation or when incentives  were not strong or effective  enough in the past, requiring institutions to    spend more budget  on providing incentive.  Moreover, a key factor of   behavioral update, the intensity of selection \citep{key:Sigmund_selfishnes}---which determines how strongly an individual bases her decision to copy another  individual's strategy on fitness  difference---might strongly impact an institutional incentives strategy  and its cost efficiency. For instance, when selection is weak such that behavioral update is close to a random process (i.e. an imitation decision  is independent of  how large the fitness difference is), providing incentives would make little difference to cause behavioral change, however strong  it is. When selection is strong, incentives that ensure a minimum fitness advantage to cooperators would already ensure positive behavioral change.


In a  stochastic, finite population context, so far this problem has been  investigated primarily  based on  agent-based and numerical simulations \citep{chen2015first,sasaki2012take,Hanijcai2018,han2018cost,cimpeanu2019exogenous}.
Results demonstrate several  interesting phenomena, such as the significant influence of the intensity of selection on incentive strategies and optimal costs. However, there is  no satisfactory rigorous analysis available at present that allows one to derive  the optimal way of providing incentives  (for a given   cooperation dilemma and desired level of cooperation). This is a challenging problem because of the large but finite population size and the complexity of stochastic processes on the population. 

In this paper, we provide exactly such a rigorous analysis. We study cooperation dilemmas in both pairwise  (the donation game) and multi-player  (the Public Goods game) settings.  
These  are the most studied  models for studying cooperation where individually defection is always preferred over cooperation while mutual cooperation is the preferred collective outcome for the population as a whole.  
Adopting a popular stochastic evolutionary game approach for analysing well-mixed finite populations \citep{key:novaknature2004,key:imhof2005,nowak:2006bo},  we derive the total expected costs of providing institutional reward or punishment, characterize their asymptotic limits (infinite population, weak selection and large selection) and show the existence of a phase transition phenomena in the optimization problem when the intensity of  selection varies. We calculate the critical threshold of phase transitions and study the minimization problem when the intensity of selection is under and above the critical value. We furthermore provide numerical simulations to demonstrate the analytical results.

The rest of the paper is organized as follows. In Section \ref{sec: models} we introduce the models and methods deriving mathematical optimization problems that will be studied. The main results of the paper are presented in Section \ref{sec: main result}. In Section \ref{sec: numerics} we give explicit computations for small populations and numerical investigations demonstrating the analytical results. In Section \ref{sec: summary} we discuss possible extensions for future work. Finally, detailed computations, technical lemmas and proofs the main results are provided in the Supplementary Information attached at the end of the paper.

\section{Materials and Methods} 
\label{sec: models}
\subsection{Cooperation dilemmas} 
We consider a well-mixed, finite population of $N$   self-regarding individuals or players, who interact with each other using one of the following cooperation dilemmas, namely the Donation Game (DG) and the Public Goods Game (PGG). Let $\Pi_C(i)$ and $\Pi_D(i)$ denote the payoffs of Cooperator (C) and Defector (D) in a population with $i$ C players and $N-i$ D players. We show below that $\Pi_C(i) - \Pi_D(i)$ is a constant, which is denoted by $\delta$. For cooperation dilemmas, it's always the case that $\delta < 0$. 
\subsubsection{Donation Game (DG)} 
 In a DG,  two players decide simultaneously whether to cooperate (C) or defect (D). The payoff matrix of the DG (for row player) is given as follows  
\[
 \bordermatrix{~ & C & D\cr
                  C & b-c & -c \cr
                  D & b & 0  \cr
                 }, 
\]
where $c$ and $b$ represent the cost and benefit of cooperation, where $b > c$. DG is a special version of the Prisoner's Dilemma game (PD). 

We obtain 
\begin{equation*} 
\begin{split} 
\Pi_C(i) &=\frac{(i-1)\pi_{C,C} + (N-i)\pi_{C,D}}{N-1} = \frac{(i-1) (b-c) + (N-i) (-c)}{N-1}  ,\\
\Pi_D(i) &=\frac{i\pi_{D,C} + (N-i-1)\pi_{D,D}}{N-1} =\frac{i b}{N-1}.
\end{split}
\end{equation*} 
Thus, 
$$\delta = \Pi_C(i) - \Pi_D(i) =  -(c + \frac{b}{N-1})$$

\subsubsection{Public Goods Game (PGG)} 
In a PGG, Players interact in a group of size $n$, where they decide whether to contribute a cost $c$ to the public goods. Cooperator contributes and defector does not. The total contribution will be multiplied by a factor $r$ ($1 < r < n$), and shared equally among all members of the group.

We obtain \citep{hauert2007via}
\begin{equation*} 
\begin{split} 
\Pi_C(i) &= \sum^{n-1}_{j=0}\frac{\dbinom{i-1}{j}\dbinom{N-i}{n-1-j}}{
 \dbinom{N-1}{n-1}} \ \left(\frac{(j+1)rc}{n} - c\right) =\frac{rc}{n}\left(1 + (i-1)\frac{n-1}{N-1}\right) - c ,\\
\Pi_D(i) &=\sum^{n-1}_{j=0}\frac{\dbinom{i}{j}\dbinom{N-1-i}{n-1-j}}{
 \dbinom{N-1}{n-1}} \ \frac{jrc}{n} =\frac{rc(n-1)}{n(N-1)}i.
\end{split}
\end{equation*} 
Thus, 
$$\delta = \Pi_C(i) - \Pi_D(i) = -c \left(1 - \frac{r(N-n)}{n(N-1)} \right).
$$
\subsection{Institutional reward and punishment} 
To reward a cooperator (respectively, punish a defector), the institution has to pay an amount $\theta/a$ (resp., $\theta/b$) so that the cooperator's (defector's)  payoff  increases (decreases) by $\theta$, where $a, b > 0$ are constants representing the efficiency ratios of providing incentive.  Without losing generality (as we study reward and punishment separately), we set $a = b = 1$. 
\subsubsection{Cost of incentives}
\label{sec: cost of incentive}
We adopt here the finite population dynamics with the Fermi strategy update rule \citep{traulsen2006} (see Methods), stating   that a player $A$ with fitness $f_A$ adopts the strategy of another player $B$ with fitness $f_B$ with a probability given by, $P_{A,B}=\left(1 + e^{-\beta(f_B-f_A)}\right)^{-1}$, where    $\beta$ represents  the intensity of selection. 
We compute the expected number of times the population contains  $i$ C players, $1 \leq i \leq N-1$. For that, we consider an absorbing Markov chain of $(N+1)$ states, $\{S_0, ..., S_N\}$, where $S_i$ represents a population with $i$ C players.  $S_0$ and $S_N$ are absorbing states. Let  $U = \{u_{ij}\}_{i,j = 1}^{N-1}$ denote the transition matrix between the $N-1$ transient states, $\{S_1, ..., S_{N-1}\}$. The transition probabilities can be defined as follows, for $1\leq i \leq N-1$: 
\begin{equation} 
\begin{split} 
u_{i,i\pm j} &= 0 \qquad \text{ for all } j \geq 2, \\
u_{i,i\pm1} &= \frac{N-i}{N} \frac{i}{N} \left(1 + e^{\mp\beta[\Pi_C(i) - \Pi_D(i)+\theta]}\right)^{-1},\\
u_{i,i} &= 1 - u_{i,i+1} -u_{i,i-1},
\end{split} 
\end{equation}
 where $\Pi_C(i)$ and $\Pi_D(i)$ represent the average payoffs of a C and D player, respectively, in a population with $i$ C players and $(N-i)$ D players. 
 
The entries $n_{ij}$ of the so-called fundamental matrix $\mathcal{N}=(n_{ij})_{i,j=1}^{N-1}= (I-U)^{-1}$ of the absorbing Markov chain gives the expected number of times the population is in the state $S_j$ if it is started in the transient state $S_i$ \citep{kemeny1976finite}.
As a mutant can  randomly occur either at $S_0$ or $S_N$, the expected number of visits at state $S_i$ is: $\frac{1}{2} (n_{1i} + n_{N-1,i})$.

The total cost per generation is
\begin{equation*}
\theta_i = \begin{cases}
i \times \theta\quad\text{in the case of institutional reward},\\
(N-i) \times \theta\quad\text{in the case of institutional punishment}.
\end{cases}
\end{equation*}
Hence, the expected total cost of interference for    institutional reward and institutional punishment are respectively 
\begin{equation} 
\label{eq:total_investment}
E_r(\theta) = 
\frac{\theta}{2} \sum_{i=1}^{N-1}(n_{1i} + n_{N-1,i}) i,\quad \text{and}\quad E_p(\theta)=\frac{\theta}{2} \sum_{i=1}^{N-1}(n_{1i} + n_{N-1,i}) (N-i).
\end{equation}
\subsubsection{Cooperation frequency}
Since the population consists of only two strategies, the fixation  probabilities of a C (respectively, D) player in a (homogeneous) population of D (respectively, C) players  when the interference scheme is carried out are (see Methods), respectively, 
\begin{equation} 
\label{eq:fixprob} 
\begin{split}
\rho_{D,C} &= \left(1+\sum_{i = 1}^{N-1} \prod_{k = 1}^i \frac{1+e^{\beta(\Pi_k(C)-\Pi_k(D) + \theta)}}{1+e^{-\beta(\Pi_k(C)-\Pi_k(D)+\theta)}}  \right)^{-1}, \\
\rho_{C,D} &= \left(1+\sum_{i = 1}^{N-1} \prod_{k = 1}^i \frac{1+e^{\beta(\Pi_k(D)-\Pi_k(C) - \theta)}}{1+e^{-\beta(\Pi_k(D)-\Pi_k(C)-\theta)}}  \right)^{-1}.
\end{split}
\end{equation} 
Computing the stationary distribution using  these fixation probabilities, we  obtain    the frequency of cooperation (see Methods) $$\frac{\rho_{D,C}}{\rho_{D,C}+\rho_{C,D}}.$$
Hence, this frequency of cooperation can be maximised by maximising 
\begin{equation}
\label{eq:max}
\max_{\theta} \left(\rho_{D,C}/\rho_{C,D}\right).  
\end{equation} 
 The fraction in Equation ~\eqref{eq:max} can be simplified as follows \citep{nowak:2006bo} 
\begin{eqnarray}
\nonumber
\frac{\rho_{D,C}}{\rho_{C,D}} &=&  \prod_{k = 1}^{N-1} \frac{T^-(k)}{T^+(k)} =\prod_{k = 1}^{N-1} \frac{1 + e^{\beta[\Pi_k(C)-\Pi_k(D) + \theta]}}{1 + e^{-\beta[\Pi_k(C)-\Pi_k(D) + \theta]}} \\
\nonumber
&=& e^{\beta\sum_{k = 1}^{N-1} \left(\Pi_k(C)-\Pi_k(D) + \theta\right)} \\
\label{eq:max_Q_prime}
 &=& e^{\beta (N-1)(\delta +  \theta)}. 
 \end{eqnarray}
 In the above transformation, $T^-(k)$ and $T^+(k)$ are the probabilities  to increase or decrease the number  of C players  (i.e. $k$) by one in each time step, respectively (see Methods for details). 

We consider non-neutral selection, i.e.  $\beta > 0$ (under neutral selection, there is no need to use incentives). Assuming that we desire to obtain  at least an $\omega  \in [0,1]$ fraction of cooperation, i.e. $\frac{\rho_{D,C}}{\rho_{D,C}+\rho_{C,D}} \geq \omega$, it follows from Equation ~\eqref{eq:max_Q_prime}  that
\begin{equation} 
\label{eq:omega_fraction}
 \theta \geq \theta_0(\omega) = \frac{1}{(N-1)\beta} \log\left(\frac{\omega}{1-\omega}\right) - \delta .
\end{equation}
Therefore it is guaranteed that if $\theta  \geq \theta_0(\omega)$,  at least an $\omega$ fraction of cooperation can be expected.
From this  condition it implies  that the lower bound of $\theta$ monotonically  depends on $\beta$. Namely,  when $\omega \geq 0.5$,  it increases with $\beta$ while decreases for $\omega < 0.5$. 

\subsubsection{Optimization problems} 
Bringing all ingredients together, we obtain the following cost-optimization problems of institutional incentives in stochastic finite populations
\begin{equation}
\label{eq: optimisation}
\min_{\theta\geq\theta_0(\omega)}E(\theta),
\end{equation}
where $E$ is either $E_r$ or $E_p$, defined in \eqref{eq:total_investment}, which respectively corresponds to institutional reward and punishment.

\subsection{Methods: Evolutionary Dynamics in Finite Populations} 
Both the analytical and numerical results obtained here use  Evolutionary Game Theory (EGT) methods for finite populations \citep{key:novaknature2004,key:imhof2005,nowak:2006bo}.  In such a setting,  players'  payoff represents their \emph{fitness} or social \emph{success}, and  evolutionary dynamics is shaped  by social learning \citep{key:Hofbauer1998,key:Sigmund_selfishnes}, whereby the  most successful players will tend to be imitated more often by the other players. In the current work, social learning is modeled using  the so-called pairwise comparison rule \citep{traulsen2006},  assuming  that a player $A$ with fitness $f_A$ adopts the strategy of another player $B$ with fitness $f_B$ with probability given by the Fermi function, 
$P_{A,B}=\left(1 + e^{-\beta(f_B-f_A)}\right)^{-1}$,  
where  $\beta$ conveniently describes the selection intensity ($\beta=0$ represents neutral drift while $\beta \rightarrow \infty$ represents increasingly deterministic selection). 
 
 For convenience of numerical computations, but without affecting analytical results, we assume here small mutation limit\citep{Fudenberg2005,key:imhof2005,key:novaknature2004}. 
As such,  at most two strategies are present in the population simultaneously, and  the behavioural dynamics can thus be  described by a Markov Chain, where each state represents a homogeneous population and  the transition probabilities between any two states  are given by the fixation probability of a single mutant \citep{Fudenberg2005,key:imhof2005,key:novaknature2004}. The resulting Markov Chain has a stationary distribution, which describes the average time the population spends in an  end state.

The fixation probability that a single mutant  A taking over a whole population with  $(N-1)$ B players is as follows (see e.g. references for details   \citep{traulsen2006,Karlin:book:1975, key:novaknature2004})
\begin{equation} 
\label{eq:fixprob} 
\rho_{B,A} = \left(1 + \sum_{i = 1}^{N-1} \prod_{j = 1}^i \frac{T^-(j)}{T^+(j)}\right)^{-1},
\end{equation} 
where $T^{\pm}(k) =  \frac{N-k}{N} \frac{k}{N} \left[1 + e^{\mp\beta[\Pi_A(k) - \Pi_B(k)]}\right]^{-1}$ describes  the probability to change the number of A players  by $\pm$ one in a time step.
Specifically, when $\beta = 0$, $\rho_{B,A} = 1/N$, representing the transition probability at neural limit. 

Having obtained the fixation probabilities between any two states of a Markov chain, we  can  now describe its stationary distribution. Namely, considering a set of $s$ strategies,  $\{1,...,s\}$, their stationary distribution  is given by the normalised eigenvector associated with the eigenvalue $1$ of the transposed of  a matrix $M = \{T_{ij}\}_{i,j = 1}^s$, where $T_{ij, j \neq i} = \rho_{ji}/(s-1)$ and  $T_{ii} = 1 - \sum^{s}_{j = 1, j \neq i} T_{ij}$.  (See e.g. references \citep{Fudenberg2005,key:imhof2005} for further details). 

\section{Main results}
\label{sec: main result}
The present paper provides a rigorous analysis for the expected total cost of providing institutional incentive \eqref{eq:total_investment} and the associated optimization problem \eqref{eq: optimisation}. In the following theorems, $E$ denotes the cost function either for  institutional reward, $E_r$, or  institutional punishment, $E_p$, as  obtained in \eqref{eq:total_investment}. Also,  $H_N$ denotes the well-known harmonic number
\begin{equation}
\label{eq: harmonic number}
H_N:=\sum_{j=1}^{N-1}\frac{1}{j}.
\end{equation}
Our first main result provides qualitative properties and asymptotic limits of $E$.
\begin{theorem}[Qualitative properties and  asymptotic limits of  expected total cost functions]
\label{thm: main theorem 1}
\begin{enumerate}[(1)]\
\item (regularity) $\theta\mapsto E(\theta)$ is a smooth function on $\mathbb{R}$.
\item (finite population estimates) The expected total cost of providing incentive satisfies the following estimates for all finite population of size $N$
\begin{equation}
\label{eq: estimates}
\frac{N^2\theta}{2}\Big(H_N+\frac{1}{N-1}\Big)\leq E(\theta)\leq N(N-1)\theta \Big(H_N+1\Big).
\end{equation}
\item (infinite population limit) The expected total cost of providing incentive  satisfies the following asymptotic behaviour when the population size $N$ tends to $+\infty$
\begin{equation}
\lim_{N\rightarrow +\infty}\frac{E(\theta)}{\frac{N^2\theta}{2}(\ln N+\gamma)}=\begin{cases}
1+e^{-\beta|\theta-c|} \quad \text{for DG},\\
1+e^{-\beta|\theta-c|}e^{\beta c\frac{r}{n}}\quad\text{for PGG},
\end{cases}
\end{equation}
where $\gamma=0.5772...$ is the Euler–Mascheroni constant.
\item (weak selection limits) The expected total cost of providing incentive satisfies the following asymptotic limit when the selection strength $\beta$ tends to $0$:
\begin{equation}
\label{eq: zero-intensity limit punishment}
\lim\limits_{\beta\rightarrow 0}E(\theta)=N^2 \theta H_N.
\end{equation}
\item (large selection limits) The expected total cost of providing incentive satisfies the following asymptotic limit when the selection strength $\beta$  tends to $+\infty$
\begin{align}
\label{eq: large selection limits}
&\lim_{\beta\rightarrow +\infty}E_r(\theta)=\begin{cases}
\frac{N^2}{2}\theta\Big(\frac{1}{N-1}+H_N\Big)\quad \text{for}\quad \theta<-\delta,\\
N^2\theta H_N \quad \text{for}\quad \theta=-\delta,\\
\frac{N^2}{2}\theta \Big(1+H_N) \quad \text{for}\quad \theta>-\delta.
\end{cases}\\
&\lim\limits_{\beta\rightarrow+\infty}E_p(\theta)=\begin{cases}
\frac{N^2\theta}{2}\Big(1+H_N\Big) \quad\text{for}\quad \theta<-\delta,\\
N^2\theta H_N \quad \text{for}\quad \theta=-\delta,\\
\frac{N^2\theta}{2}\Big(H_N+\frac{1}{N-1}\Big)\quad \text{for}\quad \theta>-\delta.
\end{cases}
\end{align}
\item (difference between reward and punishment costs) The difference between the expected total costs of reward and punishment  is given by
\begin{equation}
\label{eq: difference between costs}
(E_r-E_p)(\theta)=\begin{cases}
<0,\quad\text{for}\quad \theta<-\delta,\\
=0,\quad\text{for}\quad \theta=-\delta,\\
>0,\quad \text{for}\quad \theta>-\delta.
\end{cases}
\end{equation}
\end{enumerate}
\end{theorem}
Our second main result concerns the optimization problem \eqref{eq: optimisation}. We show that the cost function $E$ exhibits a phase transition when the strength of the selection $\beta$ varies.
\begin{theorem}[Optimization problems and phase transition phenomenon]
\label{thm: main theorem}
\begin{enumerate}[(1)]\
\item(phase transition phenomena and behaviour under the threshold) Define
$$
F^*=\begin{cases}
\min\{F(u): P(u)>0\}\quad\text{in the reward case},\\
\min\{\hat{F}(u): \hat{P}(u)>0\}\quad\text{in the punishment case,}
\end{cases}
$$
where $P(u)$, $F(u)$, $\hat{P}$ and $\hat{F}$ can be explicitly found in the Supporting Information. There exists a threshold value $\beta^*$ given by 
$$
\beta^*=-\frac{F^*}{\delta}>0,
$$
such that $\theta\mapsto E(\theta)$ is non-decreasing for all $\beta\leq \beta^*$ and is non-monotonic when $\beta>\beta^*$. As a consequence, for $\beta\leq \beta^*$ 
\begin{equation}
\label{eq: smallbeta}
\min\limits_{\theta\geq\theta_0}E(\theta)=E(\theta_0).
\end{equation}  
\item (behaviour above the threshold value) For $\beta>\beta^*$, the number of changes of the sign of $E'(\theta)$ is at least two for all $N$ and is exactly two for $N\leq N_0=100$. As a consequence,  
for $N\leq N_0$, there exist $\theta_1<\theta_2$ such that for $\beta>\beta^*$, $E(\theta)$ is increasing when $\theta< \theta_1$, decreasing when $\theta_1<\theta<\theta$ and increasing when $\theta>\theta_2$. Thus, for $N\leq N_0$, 
$$
\min\limits_{\theta\geq\theta_0}E(\theta)=\min\{E(\theta_0),E(\theta_2)\}.
$$
\end{enumerate}
\end{theorem}
Based on numerical simulations, we conjecture that the requirement that $N\leq N_0$ could be removed and Theorem \ref{thm: main theorem} is true for all finite $N$. Theorem \ref{thm: main theorem} gives rise to the following algorithm to determine the optimal value $\theta^*$ for $N\leq N_0$.
\begin{algorithm}[\textbf{Finding optimal cost of incentive $\theta^\star$}]\ \\ \ \\
\label{algo}
\textbf{Inputs}: i) $N\leq N_0$: population size, ii) $\beta$: intensity of selection, iii) game and parameters: PD ($c$ and $b$) or PGG ($c$, $r$ and $n$), iv) $\omega$: minimum desired  cooperation level. 
\begin{enumerate}[(1)]
\item Compute $\delta$ \Big\{in PD: $\delta = - (c + \frac{b}{N-1})$; in PGG: $\delta = -c \left(1 - \frac{r(N-n)}{n(N-1)} \right)$\Big\}.
\item Compute $\theta_0 = \frac{1}{(N-1)\beta} \log\left(\frac{\omega}{1-\omega}\right) - \delta$; 
\item Compute 
$$
F^*=\begin{cases}
\min\{F(u): P(u)>0\}\quad\text{in the reward case},\\
\min\{\hat{F}(u): \hat{P}(u)>0\}\quad\text{in the punishment case,}
\end{cases}
$$
where $P(u)$, $F(u)$, $\hat{P}(u)$ and $\hat{F}(u)$ are explicitly defined in the Supporting Information.
\item Compute $\beta^*=-\frac{F^*}{\delta}$.
\item If $\beta\leq \beta^*$:
$$
\theta^*=\theta_0,\quad \min E(\theta)=E(\theta_0).
$$ 
\item Otherwise (i.e. if $\beta>\beta^*$)
\begin{enumerate}
\item Compute $u_2$ that is the largest root of the equation $F(u)+\beta \delta=0$ for the reward case or that of $\hat{F}(u)+\beta \delta=0$ for the punishment case. 
\item Compute $\theta_2=\frac{\log u_2}{\beta}-\delta$.
	\begin{itemize} 
		\item If $\theta_2 \leq \theta_0$: $\theta^*=\theta_0,\quad \min E(\theta)=E(\theta_0)$;
		\item Otherwise (if  $\theta_2 > \theta_0$): 
			\begin{itemize}  
				\item If $E(\theta_0)\leq E(\theta_2)$: $\theta^*=\theta_0,\quad \min E(\theta)=E(\theta_0)$;
				\item if $E(\theta_2)< E(\theta_0)$: $\theta^*=\theta_2,\quad \min E(\theta)=E(\theta_2)$.
			\end{itemize} 
	\end{itemize} 
\end{enumerate}
\end{enumerate}
\noindent \textbf{Output}: $\theta^*$ and $E(\theta^*)$.
\end{algorithm} 
\section{Explicit computations for some small $N$ and numerical investigations}
\label{sec: numerics}
In this section, we will present explicit computations for some small $N$ and provide several numerical investigations to demonstrate the main analytical results, Theorem \ref{thm: main theorem} and Algorithm \ref{algo}.
\subsection{Illustrative examples}
\begin{example}
Let us consider the case $N=3$, the total cost of interference for institutional reward and institutional punishment are respectively given by
\begin{align*}
E_r(\theta)=\frac{9\theta}{4}\Big(5+\frac{4e^x-1}{1+e^x+e^{2x}}\Big),\quad E_p(\theta)=\frac{9\theta}{4}\Big(4+\frac{5e^x+1}{1+e^x+e^{2x}}\Big),
\end{align*}
where $x=\beta(\theta+\delta)$. For institutional reward, 
$$
F(u)=\frac{(u+1)(5u+4)(u^2+u+1)}{u(2u^2-u-5/2)}-\log(u),\quad P(u)=4u^2-2u-5.
$$
For institutional punishment
$$
\hat{F}(u)=\frac{(u+1)(4u+5)(u^2+u+1)}{u (\frac{5}{2}u^2+u-2)}-\log(u),\quad \hat{P}(u)=\frac{5}{2}u^2+u-2.
$$
To illustrate the phase transition phenomena, we consider the institutional reward (institutional punishment is similar). We compute $F^*$: 
$$
F^*=\min\{F(u): P(u)>0\}=10.9291 \quad\text{obtained at}\quad u^*=4.29712.
$$
Moreover, $F(u)$ is decreasing when $u<u^*$ and is increasing when $u>u^*$. The number of solutions of the equation $F(u)=\beta a$ and the sign of  $E_r'(\theta)$ depends on the value of $\beta a$.
\begin{enumerate}[(i)]
\item If $\beta a<F^*$, then $F(u)=\beta a$ has no solution and $E_r'(\theta)>0$. Hence $E_r$ is increasing.
\item If $\beta a= F^*$, then $F(u)=\beta a$ has a unique solution, and $E'(\theta)\geq 0$ and $E_r$ is non-decreasing.
\item If $\beta a>F^*$, then $F(u)=\beta a$ has two solutions $(1+\sqrt{21})/4<u_1<u_2$. 
$$
E_r'(\theta)\begin{cases}
>0\qquad u<u_1,\\
<0\qquad u_1<u<u_2,\\
>0\qquad u>u_2.
\end{cases}
$$
\end{enumerate}
Thus $E_r$ is increasing when $u<u_1$, is decreasing when $u_1<u<u_2$ and is increasing again when $u>u_2$. 

\paragraph{Algorithm for $N=3$}
\begin{enumerate}[1.]
\item Compute 
\begin{equation}
\delta=\begin{cases}
-(c+\frac{b}{N-1})\quad\text{PD games},\\
-c\Big(1-\frac{r(N-n)}{n(N-1)}\Big)\quad\text{PGG games}
\end{cases}
\end{equation}
\item Compute $\theta_0 = \frac{1}{(N-1)\beta} \log\left(\frac{\omega}{1-\omega}\right) - \delta$;
\item Compute $\beta^*=-\frac{f^*}{\delta}=-\frac{10.9291}{\delta}$.
\item If $\beta\leq \beta^*$ then
$$
\theta=\theta_0, \quad \min E(\theta)=E(\theta_0).
$$
\item If $\beta>\beta^*$
\begin{enumerate}
\item Compute $\theta_2=\frac{\log u_2}{\beta}-\delta$ where $u_2$ is the largest zero of the equation $F(u)+\beta \delta=0$.
\item If $\theta_0\geq \theta_2$: $\theta^*=\theta_0,\quad \min E(\theta)=E(\theta_0)$;
\item if $\theta_0< \theta_2$, compare $E(\theta_0)$ and $E(\theta_2)$: if $E(\theta_0)\leq E(\theta_2)$: $\theta^*=\theta_0,\quad \min E(\theta)=E(\theta_0)$; if $E(\theta_2)< E(\theta_0)$: $\theta^*=\theta_2,\quad \min E(\theta)=E(\theta_2)$.
\end{enumerate}
\end{enumerate}
\end{example}
\begin{figure}[!htb]
\centering
\includegraphics[width = \linewidth]{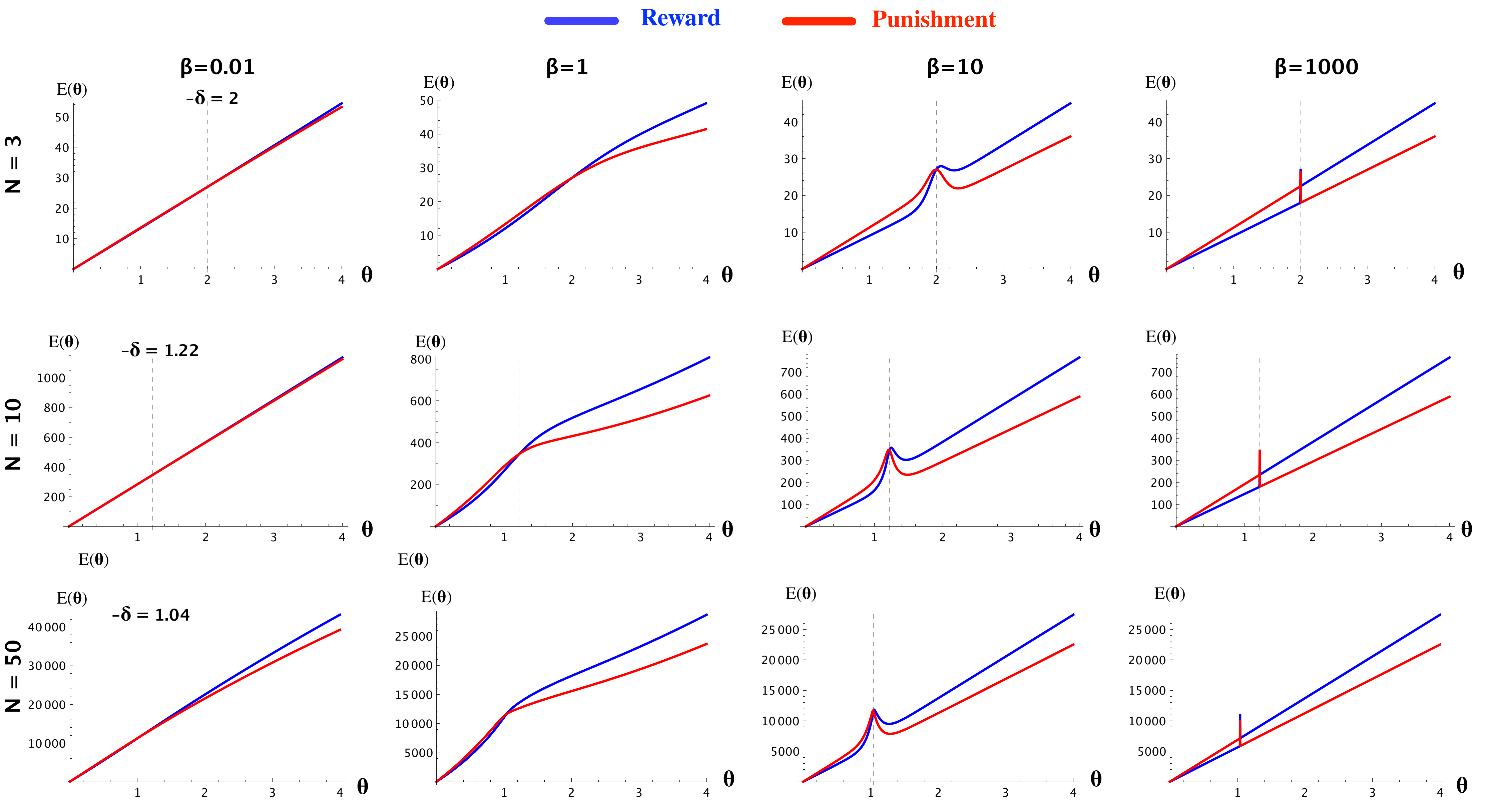}
\caption{\textbf{Compare the expected total  cost of investment $E$ for reward and punishment, for varying $\theta$}. Reward is less costly than punishment ($E_r <  E_p$)  for small $\theta$ and vice versa otherwise. The threshold of $\theta$ for this change was obtained analytically (see Theorem 1), which is exactly equal to $-\delta$.   Other parameters:  Donation Game with $b = 2$, $c = 1$.}
\label{fig:reward_vs_punishment_theta_threshold}
\end{figure} 

\begin{figure}[!htb]
\centering
\includegraphics[width = \linewidth]{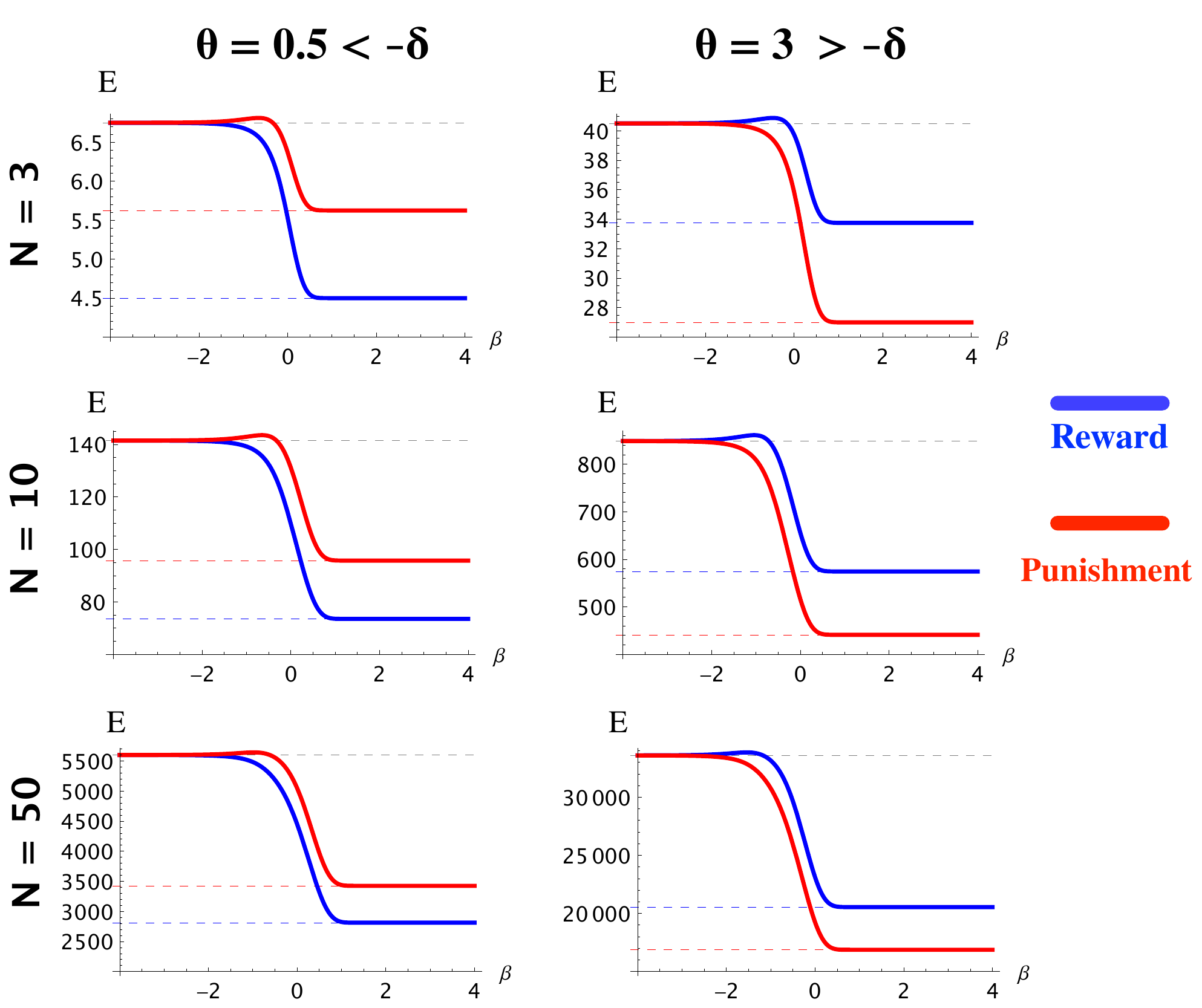}
\caption{\textbf{Total  cost of incentive $E$ for reward and punishment, for different regimes of intensity of selection, $\beta$ (Log-10 scaled), from weak to strong selection limits}. We show for different values of $N$. The dashed lines represent the corresponding theoretical limiting values obtained in Theorem \ref{thm: main theorem 1} for weak and strong selection limits. We observe that numerical results are in close accordance with those obtained theoretically.   Other parameters:  Donation Game with $b = 2$, $c = 1$.}
\label{fig:reward_vs_punishment_vary-beta}
\end{figure} 

\subsection{Numerical evaluation}
We perform several numerical investigations to demonstrate the main theoretical results. To begin with, Figure \ref{fig:reward_vs_punishment_theta_threshold} shows the expected total costs for reward and punishment (DG), for varying $\theta$. We observe that reward is less costly than punishment ($E_r <  E_p$)  for  $\theta < -\delta$ and vice versa when $\theta > -\delta$. It is exactly as shown  analytically in  Theorem \ref{thm: main theorem 1}.  This analytical result is confirmed here for different population size $N$ and intensity of selection $\beta$. 

In Figure \ref{fig:reward_vs_punishment_vary-beta} we calculate the total costs of incentive for both reward and punishment for different regimes of intensity of selection. We observe that for both weak and strong limits of selection, the theoretical results obtained in Theorem \ref{thm: main theorem 1} are confirmed, for different population sizes.

\begin{figure}[!htb]
\centering
\includegraphics[width = 0.95\linewidth]{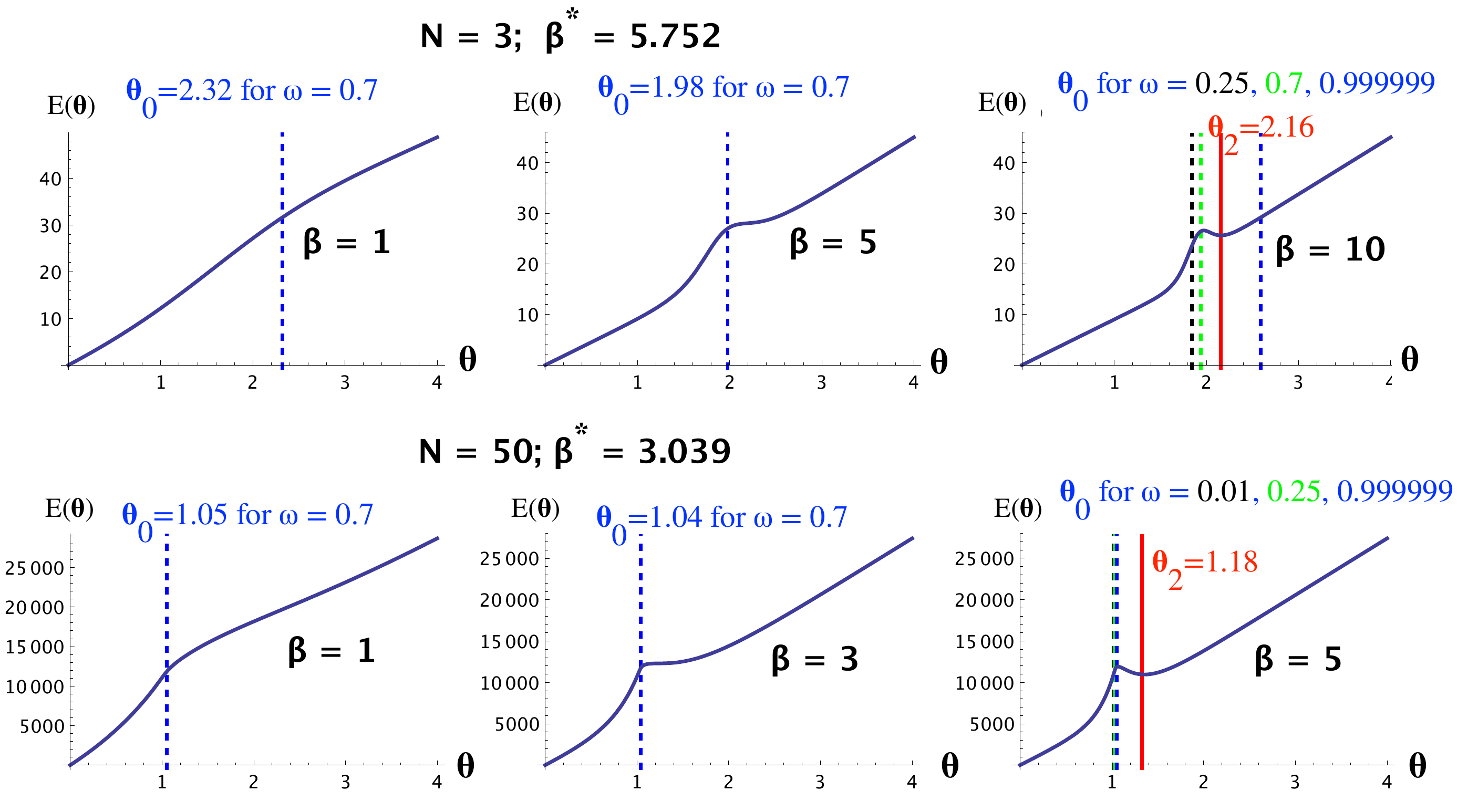}
\caption{\textbf{Using Algorithm \ref{algo} to find optimal $\theta$ that minimizes $E(\theta)$ (for institutional reward) while ensuring a minimum  level of cooperation $\omega$}. We use as examples a small population size ($N = 3$, \textbf{top row}) and a larger one ($N = 50$, \textbf{bottom row}), for  Donation Game ($b = 1.8$, $c = 1$). 
}
\label{fig:sn}
\end{figure}

 Now for Theorem \ref{thm: main theorem} and Algorithm \ref{algo}, we focus on reward for illustration.  Figure \ref{fig:sn} plots the cost function $E(\theta)$ (for institutional reward) in terms of $\theta$ for different values of $N$, $\beta$ and $\omega$ for illustrating the phase transition when varying $\beta$, in a DG. We can see that in all cases, these numerical observations are in close accordance with theoretical results. 
For example, with $N = 3$ (see top row), we found $\beta^\star =  f^\star / \delta  = 10.9291/1.9 = 5.752$. Similarly, in the second case, $\beta^\star = 3.15/1.03673  = 3.039$. For $\beta < \beta^\star$,  $E(\theta)$  are increasing functions of $\theta$. Thus, the optimal cost of incentive $\theta^\star = \theta_0$, for a given required minimum level of cooperation $\omega$. For example, with $N = 3$, for $\beta = 1$ to ensure at least 70\% of cooperation ($\omega = 0.7$), then $\theta^\star = \theta_0 = 2.32$. When $\beta \geq  \beta^\star$ one needs to compare $E(\theta_0)$ and $E(\theta_2)$. For example, with $N = 3$, $\beta = 10$: for $\omega = 0.25$  (black dashed line), then $E(\theta_0) = 23.602 < 25.6124 = EC(\theta_2)$, so  $\theta^\star = \theta_0 = 1.845$; for  $\omega = 0.7$ (green dashed line), then $E(\theta_0) = 26.446 > 25.6124 = EC(\theta_2)$, so  $\theta^\star = \theta_2 = 2.16$ (red solid line); for $\omega = 0.999999$ (blue dashed line), since $\theta_2 < \theta_0$, $\theta^\star = \theta_0 = 2.59078$.

Similarly, with a larger population size ($N = 50$), we obtained $\beta^\star = 3.039$.  In general, similar observations are obtained as in case of a small population size $N = 3$. Except that when $N$ is large, the values of $\theta_0$ for different non-extreme values of minimum required cooperation $\omega$ (say, $\omega \in (0.01, 0.99)$) is very small (given the log scale of $\omega/(1-\omega)$ in the formula of $\omega_0$). This value is also smaller than $\theta_0$, with a cost $E(\theta_0) > E(\theta_2)$, making $\theta_2$ the optimal cost of incentive.  Similar results are obtained for PGG (see Figure 2 in the Supporting Information). When $\omega$ is  extremely high  (i.e. $> 1-10^{-k}$, for a large $k$) (we don't look at extremely low value since we would like to ensure at least a sufficient  level of cooperation), then we can also see other scenarios where the optimal cost is $\theta_0$ (see Figure 1 in the SI, bottom row). We thus can observe  that for  
$\omega \in (0.01, 0.99)$, for sufficiently large population size $N$ and large enough $\beta$ ($\beta > \beta^\star + \text{a bit more}$), then the optimal value of $\omega$ is always $\theta_2$. Otherwise,   $\theta_0$ is the optimal cost.

Figure 3 in the SI plots the increase in the costs $\theta_0$, $\theta^\star$ and $E$ in order to increase the level of cooperation from $\epsilon$ to $1 - \epsilon$.

\begin{figure}[!htb]
\centering
\includegraphics[width = \linewidth]{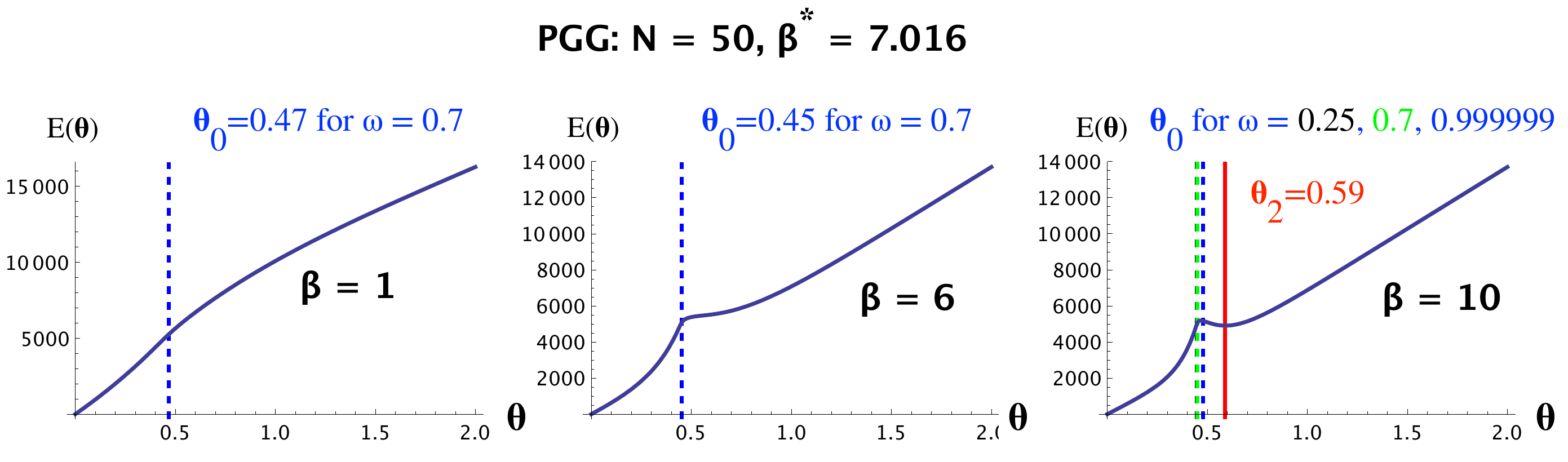}
\caption{\textbf{Using Algorithm  \ref{algo} to find optimal $\theta$ that minimizes $E(\theta)$ while ensuring a minimum  level of cooperation $\omega$, for PGG  ($r = 3$, $n = 5$, $c = 1$) with $N = 50$}. Similar observations to those with Donation Game,  are obtained.    }
\label{fig:sn-PGG}
\end{figure}

\section{Discussion}
\label{sec: summary}
Institutional incentives such as punishment and reward provide an effective tool for promoting the evolution of cooperation in social dilemmas. 
Both theoretical and experimental analysis has been provided  \citep{gurerk2006competitive,sasaki2012take,garcia2019evolution,baldassarri2011centralized,Dong2019,bardhan2005institutions,wu2014role}. 
However, past research usually ignores the question of how institutions' overall spending, i.e. the total cost of providing these incentives, can be minimized (while ensuring a desired level of cooperation). Answering this question allows one to estimate exactly how  incentives should be provided, i.e. how much to reward a cooperator and how severely to punish a wrongdoer. 
Existing works that consider this question usually omit the stochastic effects in the population, namely, the intensity of selection. 

Resorting to a stochastic evolutionary game approach for finite, well-mixed populations, we provide here theoretical results for the optimal cost of incentives that ensure a desired level of cooperation while minimizing the total budget, for a given intensity of selection, $\beta$. We show that this cost strongly depends on the value of $\beta$, due to  the existence of a phase transition in the cost functions when $\beta$ varies. This behavior is missing in works that consider a deterministic evolutionary approach \citep{wang2019exploring}.   The intensity of selection plays an important role in  evolutionary processes. Its value   differs  depending on the payoff structure (i.e., scaling game payoff matrix by a factor is equivalent  to dividing $\beta$ by that factor) and specific populations, which can be estimated  in behavioral  experiments  \citep{traulsen2010human,randUltimatum,zisisSciRep2015,domingos2020timing}. Thus, our analysis provides a way to calculate the optimal incentive cost for a given population and interaction game at hand.  

As of theoretical importance, we  characterize  asymptotic behaviors of the total cost functions for both reward and punishment (namely,  large population, weak selection and strong selection, limits) and compare these functions for the two types of incentive. We show that punishment is alway more costly for small (individual) incentive  cost ($\theta$) but less so when this cost is above a certain threshold. This result provides  insights into the choice of which type of incentives to use. We provide an exact formula for this threshold.  Moreover, we provide numerical validation of the theoretical results.

In the context of institutional incentives modelling, a crucial issue is the question of how to maintain the budget of incentives providing. The problem of who pays or contributes to the budget is a social dilemma itself, and how to escape this dilemma is critical research question. In this work we focus on the question of how to optimize the budget used for provided incentives. For future research, we plan to generalise the results of this paper to more complex scenarios such as other social dilemmas, mixed incentives, network effects, and any mutation.


\section{Acknowledgements}
T.A.H. acknowledges support from Leverhulme Research Fellowship (RF-2020-603/9) and   Future of Life Institute (grant RFP2-154). 
\newpage
\bibliographystyle{plainnat}
\bibliography{ref}


\section*{Supplementary material}

# Supplementary Information: Cost efficiency of institutional incentives in finite populations

Manh Hong Duong[1] and The Anh Han[2]

[1]School of Mathematics, University of Birmingham, Birmingham B15 2TT, UK.
Email: h.duong@bham.ac.uk

[2]School of Computing, Engineering and Digital Technologies, Teesside University, TS1 3BX, UK. Email: T.Han@tees.ac.uk

March 1, 2021

**Abstract**

In this Supplementary Information, we provide detailed calculations and proofs as well as illustrative figures for the analytical results in the main text of the present paper.

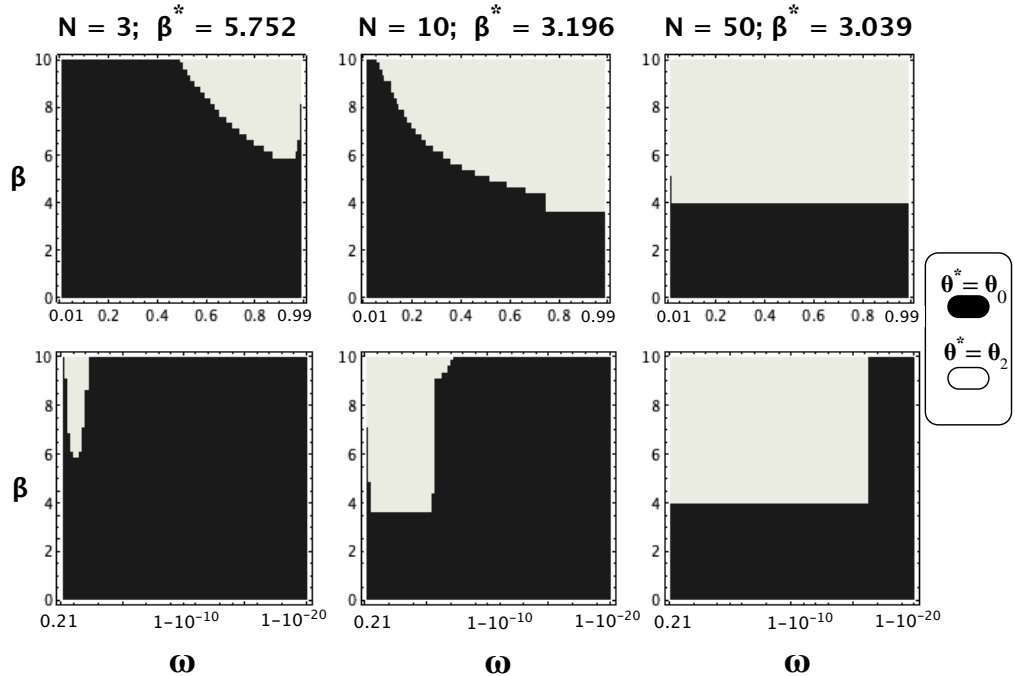

Figure 1: **Optimal value of the incentive cost $\theta^\star$, for varying the minimum required cooperation level $\omega$ and the intensity of selection $\beta$.** Left column: $N = 3$; Middle column: $N = 10$; Right column: $N = 50$. Top row shows results for a non-extreme value range of $\omega$: $\omega \in [0.01, 0.99]$. Bottom row shows results in a log scale to examine values of $\omega$ extremely close to 1. Other parameters: Donation Game ($b = 1.8$, $c = 1$).

Figure 2: **Plot of $\beta^\star$ for varying the population size** $N$. We consider three games Donation Game and Public Goods Game ($n = 5$), varying their own parameters as well (in both cases, $c = 1$). The right-hand side figures plot the value $-1/\delta$ as a function of $N$. This is the upper bound for the critical threshold $\beta^*$ in the institutional reward.

# 1 Detailed computations: reward

## 1.1 Computation of the fundamental matrix

The key in our analysis is that the fundamental matrix can be computed explicitly based on the following result.

**Lemma 1.1.** *(Inverse of a tridiagonal matrix [Huang and McColl, 1997]) Let $T$ be a general Toeplitz tridiagonal matrix of the form*

$$T = \begin{pmatrix} 1 & -u & & & & & \\ -l & 1 & -u & & & & \\ & \ddots & \ddots & \ddots & & & \\ & & -l & 1 & -u & & \\ & & & \ddots & \ddots & \ddots & \\ & & & & -l & 1 & -u \\ & & & & & -l & 1 \end{pmatrix} \tag{1}$$

2

Figure 3: **Plot of the increase in the costs $\theta_0$, $\theta^\star$ and $E$ in order to increase the minimum required level of cooperation from $\epsilon$ to $1-\epsilon$ (i.e. $1-2\epsilon$ increase) for varying the intensity of selection, $\beta$ ($N=3$, donation game with $b=1.8$, $c=1$).** In general, all three costs decrease with $\beta$. For a sufficiently weak selection (i.e. small $\beta$, see **top row**), a significant increase in the costs are required to increase even a rather small amount of cooperation (e.g., when $\epsilon = 0.4$, or an increase of 0.2 of cooperation). In contrast, for a sufficiently strong selection (see **bottom row**), see the lower row for a clearer view), even a big increase of cooperation (when $\epsilon = 0.01$, or an increase of 0.98), a small increase in the costs is required. These observations suggest that the larger $\beta$, the overall cooperation outcome changes more significantly for the same change in the per-individual investment; thereby requiring more careful in choosing this cost to ensure a low, optimal cost.

Then the inverse of $T$ is given by

$$(T^{-1})_{ij} = \frac{(\lambda_+^i - \lambda_-^i)(\lambda_+^{n-j+1} - \lambda_-^{n-j+1})}{(\lambda_+ - \lambda_-)(\lambda_+^{n+1} - \lambda_-^{n+1})} u^{j-i} \qquad i < j, \qquad (2)$$

$$(T^{-1})_{ij} = \frac{(\lambda_+^j - \lambda_-^j)(\lambda_+^{n-i+1} - \lambda_-^{n-i+1})}{(\lambda_+ - \lambda_-)(\lambda_+^{n+1} - \lambda_-^{n+1})} l^{j-i} \qquad i \geq j, \qquad (3)$$

where

$$\lambda_+ = \frac{1 + \sqrt{1-4lu}}{2}, \quad \lambda_- = \frac{1 - \sqrt{1-4lu}}{2}. \qquad (4)$$

The following lemma is a special case of Lemma 1.1 which will be used in the subsequent analysis.

**Lemma 1.2.** *Let $W$ be the following Toeplitz tri-diagonal matrix of order $N-1$*

$$W = \begin{pmatrix} 1 & -a_1 & & & & & \\ -c_1 & 1 & -a_1 & & & & \\ & \ddots & \ddots & \ddots & & & \\ & & -c_1 & 1 & -a_1 & & \\ & & & \ddots & \ddots & \ddots & \\ & & & & -c_1 & 1 & -a_1 \\ & & & & & -c_1 & 1 \end{pmatrix}, \qquad (5)$$

3

where
$$a_1 := (1 + e^{-x})^{-1}, \quad c := (1 + e^x)^{-1}, \quad x = \beta(\theta - c - \frac{b}{N-1}).$$

Then
$$(W^{-1})_{1,j} = \frac{(e^{(N-1)x} - e^{(j-1)x})(1 + e^x)}{e^{Nx} - 1} = \frac{e^{(j-1)x}(1 + e^x)(1 + e^x + \ldots + e^{(N-j-1)x})}{1 + e^x + \ldots + e^{(N-1)x}}$$
$$(W^{-1})_{N-1,j} = \frac{(e^{jx} - 1)(1 + e^x)}{e^{Nx} - 1} = \frac{(1 + e^x + \ldots + e^{(j-1)x})(1 + e^x)}{1 + e^x + \ldots + e^{(N-1)x}}.$$

*Proof.* Applying Lemma 1.1 for $u = a = (1 + e^{-x})^{-1}, l = c = (1 + e^x)^{-1}, n = N - 1$, we obtain
$$1 - 4lu = 1 - \frac{4}{(1 + e^x)(1 + e^{-x})} = \frac{e^x + e^{-x} - 2}{e^x + e^{-x} + 2} = \frac{e^{2x} - 2e^x + 1}{e^{2x} + 2e^x + 1} = \left(\frac{e^x - 1}{e^x + 1}\right)^2.$$

Hence
$$\lambda_\pm = \frac{1 \pm \sqrt{1 - 4lu}}{2} = \frac{1 \pm \left|\frac{e^x - 1}{e^x + 1}\right|}{2} = \frac{e^x + 1 \pm |e^x - 1|}{2(e^x + 1)},$$

that is
$$\lambda_+ = \frac{e^x}{e^x + 1}, \quad \lambda_- = \frac{1}{e^x + 1} \quad \text{if} \quad x \geq 0,$$
$$\lambda_+ = \frac{1}{e^x + 1}, \quad \lambda_- = \frac{e^x}{e^x + 1} \quad \text{if} \quad x < 0.$$

According to Lemma 1.1 we have
$$(W^{-1})_{1j} = \frac{\lambda_+^{n-j+1} - \lambda_-^{n-j+1}}{\lambda_+^{n+1} - \lambda_-^{n+1}} u^{j-1} = \frac{\lambda_-^{n-j+1} - \lambda_+^{n-j+1}}{\lambda_-^{n+1} - \lambda_+^{n+1}} u^{j-1},$$
$$(W^{-1})_{N-1,j} = \frac{\lambda_+^j - \lambda_-^j}{\lambda_+^{n+1} - \lambda_-^{n+1}} l^{j-1} = \frac{\lambda_-^j - \lambda_+^j}{\lambda_-^{n+1} - \lambda_+^{n+1}} l^{j-1},$$

Therefore, by swapping the role of $\lambda_+$ and $\lambda_-$ in the case $x \geq 0$ and $x < 0$, we can assume without affecting the computations that $\lambda_+ = e^x/(1 + e^x)$ and $\lambda_- = 1/(1 + e^x)$ and obtain

$$\begin{aligned}(W^{-1})_{1j} &= \frac{\left(\frac{\lambda_+}{\lambda_-}\right)^{n-j+1} - 1}{\left(\frac{\lambda_+}{\lambda_-}\right)^{n+1} - 1} \cdot \frac{\lambda_-^{n-j+1}}{\lambda_-^{n+1}} \cdot u^{j-1} \\
&= \frac{e^{(n-j+1)x} - 1}{e^{(n+1)x} - 1} \cdot \frac{(1 + e^x)^{-(n-j+1)}}{(1 + e^x)^{-(n+1)}} \cdot (1 + e^{-x})^{1-j} \\
&= \frac{e^{(n-j+1)x} - 1}{e^{(n+1)x} - 1} \cdot (1 + e^x)^j \cdot (1 + e^{-x})^{1-j} \\
&= \frac{e^{(n-j+1)x} - 1}{e^{(n+1)x} - 1} \cdot e^{jx}(1 + e^{-x})^j \cdot (1 + e^{-x})^{1-j} \\
&= \frac{e^{(n-j+1)x} - 1}{e^{(n+1)x} - 1} \cdot e^{jx}(1 + e^{-x}) \\
&= \frac{e^{(n+1)x} - e^{jx}}{e^{(n+1)x} - 1} \cdot (1 + e^{-x}) \\
&= \frac{e^{nx} - e^{(j-1)x}}{e^{(n+1)x} - 1} \cdot (1 + e^x) \\
&= \frac{e^{(N-1)x} - e^{(j-1)x}}{e^{Nx} - 1} \cdot (1 + e^x).\end{aligned}$$

Similar computations yield

$$(W^{-1})_{N-1,j} = \frac{(e^{jx} - 1)(1 + e^x)}{e^{Nx} - 1}.$$

This completes the proof of the lemma. □

## 1.2 The total cost function and its derivative

Using results from the previous section, in this section we will compute the total cost of interference and its derivative.

**Lemma 1.3.** *Recalling that $x = \beta(\theta - c - \frac{b}{N-1}) := \beta(\theta - a)$ (where $a := c + b/(N-1) = -\delta$), the total cost of interference for the institutional reward is given by*

$$E_r(\theta) = \frac{N^2 \theta(1 + e^x)}{2(1 + e^x + \ldots + e^{(N-1)x})} \left[ \left(1 + e^x + \ldots + e^{(N-2)x}\right) \sum_{j=1}^{N-1} \frac{1}{j} + e^{(N-1)x} \sum_{j=1}^{N-1} \frac{e^{-jx}}{j} \right]. \quad (6)$$

*It derivative is given by*

$$E_r'(\theta) = \frac{1}{g(x)^2} \left[ f(x)g(x) - (x + \beta a)(f(x)g'(x) - f'(x)g(x)) \right].$$

*When $(f(x)g'(x) - f'(x)g(x)) > 0$ we write $E_r'(\theta)$ as*

$$E_r'(\theta) = \frac{1}{g(x)^2}(f(x)g'(x) - f'(x)g(x))\Big(F(u) - \beta a\Big), \quad (7)$$

*where $u = e^x$, and $F$ is given in (11).*

*Proof.* Let $V := \{v_{ij}\}_{i,j=1}^{N-1} = I - U$, we have

$$\begin{aligned}
v_{i,i\pm j} &= 0 \quad \text{for all } j \geq 2, \\
v_{i,i\pm 1} &= -\frac{N-i}{N} \frac{i}{N} \left(1 + e^{\mp \beta[\delta + \theta]}\right)^{-1}, \\
v_{i,i} &= u_{i,i+1} + u_{i,i-1} = \frac{N-i}{N} \frac{i}{N},
\end{aligned} \quad (8)$$

Let $a_1 := \left(1 + e^{-\beta(\theta + \delta)}\right)^{-1}$ and $c_1 := \left(1 + e^{\beta(\theta + \delta)}\right)^{-1}$. The fundamental matrix $\mathbb{N} = (I - U)^{-1} = V^{-1}$ (see Section the Introduction in the main text) can be expressed as

$$\mathbb{N}^{-1} = W^{-1} \operatorname{diag}\left\{\frac{N^2}{(N-1)}, \frac{N^2}{2(N-2)}, \ldots, \frac{N^2}{(N-1)}\right\},$$

where $W$ is the Toeplitz tri-diagonal matrix given in (5). Therefore, from Lemma 1.2 the expected

total cost of interference for institutional reward can be expressed as

$$\begin{aligned}
E_r(\theta) &= \frac{\theta}{2} \sum_{i=1}^{N-1} (n_{1i} + n_{N-1,i})i \\
&= \frac{1}{2} \sum_{j=1}^{N-1} \frac{N^2}{(N-j)j} \Big((W^{-1})_{1j} + (W^{-1})_{N-1,j}\Big) j\theta \\
&= \frac{N^2 \theta(1+e^x)}{2(e^{Nx}-1)} \sum_{j=1}^{N-1} \frac{e^{(N-1)x} - e^{(j-1)x} + e^{jx} - 1}{N-j} \\
&= \frac{N^2 \theta(1+e^x)}{2(e^{Nx}-1)} \Bigg[ (e^{(N-1)x} - 1) \sum_{j=1}^{N-1} \frac{1}{j} + (e^x - 1) \sum_{j=1}^{N-1} \frac{e^{(j-1)x}}{N-j} \Bigg] \\
&= \frac{N^2 \theta(1+e^x)}{2(e^{Nx}-1)} \Bigg[ (e^{(N-1)x} - 1) \sum_{j=1}^{N-1} \frac{1}{j} + (e^x - 1) e^{(N-1)x} \sum_{j=1}^{N-1} \frac{e^{-jx}}{j} \Bigg] \\
&= \frac{N^2 \theta(1+e^x)}{2(1 + e^x + \ldots + e^{(N-1)x})} \Bigg[ (1 + e^x + \ldots + e^{(N-2)x}) \sum_{j=1}^{N-1} \frac{1}{j} + e^{(N-1)x} \sum_{j=1}^{N-1} \frac{e^{-jx}}{j} \Bigg] \\
&:= \frac{N^2 \theta}{2} \frac{f(x)}{g(x)},
\end{aligned}$$

where (recalling the definition of the harmonic number $H_N$)

$$f(x) = (1+e^x) \Big[ (1 + e^x + \ldots + e^{(N-2)x}) H_N + e^{(N-1)x} \sum_{j=1}^{N-1} e^{-jx}/j \Big],$$

$$g(x) = 1 + e^x + \ldots + e^{(N-1)x},$$

Since $x = \beta(\theta - a)$, we have

$$\begin{aligned}
\frac{d}{d\theta} E_r(\theta) = \beta \frac{d}{dx} E_r(\theta) &= \frac{\beta N^2}{2} \left\{ \frac{1}{\beta} \frac{f(x)}{g(x)} + \theta \Big( \frac{f'(x)}{g(x)} - \frac{f(x) g'(x)}{g(x)^2} \Big) \right\} \\
&= \frac{N^2}{2} \left\{ \left[ \frac{f(x)}{g(x)} + (x + \beta a) \Big( \frac{f'(x)}{g(x)} - \frac{f(x) g'(x)}{g(x)^2} \Big) \right] \right\} \\
&= \frac{N^2}{2} \frac{1}{g(x)^2} \Big[ f(x) g(x) - (x + \beta a)(f(x) g'(x) - f'(x) g(x)) \Big]. \quad (9)
\end{aligned}$$

Let $u := e^x$. From the formula of $f$ and $g$, it is more convenient to express the right-hand side of

6

(9) in terms of $u$. We have

$$f(x) = (1+e^x)\left[(1+e^x+\ldots+e^{(N-2)x})H_N + e^{(N-1)x}\sum_{j=1}^{N-1}e^{-jx}/j\right] = (1+u)\sum_{j=0}^{N-2}(H_N + \frac{1}{N-1-j})u^j,$$

$$f'(x) = u\sum_{j=0}^{N-2}(H_N + \frac{1}{N-1-j})u^j + (1+u)\left(H_N\sum_{j=0}^{N-2}je^{jx} + \sum_{j=0}^{N-2}\frac{j}{N-1-j}e^{jx}\right)$$

$$= \sum_{j=0}^{N-2}\left(H_N + \frac{1}{N-1-j}\right)u^j(j(1+u)+u),$$

$$g(x) = 1 + e^x + \ldots + e^{(N-1)x} = \sum_{j=0}^{N-1}u^j,$$

$$g'(x) = \sum_{j=1}^{N-1}je^{jx} = \sum_{j=1}^{N-1}ju^j.$$

Therefore

$$f(x)g'(x) - f'(x)g(x)$$

$$= (1+u)\left(\sum_{j=0}^{N-2}(H_N + \frac{1}{N-1-j})u^j\right)\left(\sum_{j=1}^{N-1}ju^j\right)$$

$$- \left(\sum_{j=0}^{N-2}\left(H_N + \frac{1}{N-1-j}\right)u^j(j(1+u)+u)\right)\left(\sum_{j=0}^{N-1}u^j\right)$$

$$= (1+u)u\left[\left(\sum_{j=0}^{N-2}(H_N + \frac{1}{N-1-j})u^j\right)\left(\sum_{j=1}^{N-1}ju^{j-1}\right) - \left(\sum_{j=1}^{N-2}\left(H_N + \frac{1}{N-1-j}\right)ju^{j-1}\right)\left(\sum_{j=0}^{N-1}u^j\right)\right]$$

$$- u\left(\sum_{j=0}^{N-2}\left(H_N + \frac{1}{N-1-j}\right)u^j\right)\left(\sum_{j=0}^{N-1}u^j\right)$$

$$:= uP(u),$$

where

$$P(u) := (1+u)\left[\left(\sum_{j=0}^{N-2}(H_N + \frac{1}{N-1-j})u^j\right)\left(\sum_{j=1}^{N-1}ju^{j-1}\right) - \left(\sum_{j=1}^{N-2}\left(H_N + \frac{1}{N-1-j}\right)ju^{j-1}\right)\left(\sum_{j=0}^{N-1}u^j\right)\right]$$

$$- \left(\sum_{j=0}^{N-2}\left(H_N + \frac{1}{N-1-j}\right)u^j\right)\left(\sum_{j=0}^{N-1}u^j\right). \qquad (10)$$

We also have

$$f(x)g(x) = (1+u)\left(\sum_{j=0}^{N-2}(H_N + \frac{1}{N-1-j})u^j\right)\left(\sum_{j=0}^{N-1}u^j\right) := Q(u).$$

Note that $x + \beta a = \beta\theta > 0$ and $f(x)g(x) > 0$. It follows that if $f(x)g'(x) - f'(x)g(x) \le 0$ then $E'_r(\theta) > 0$. Consider the case $f(x)g'(x) - f'(x)g(x) > 0$. Let $u = e^x$ and define

$$F(u) := \frac{Q(u)}{uP(u)} - \log(u). \qquad (11)$$

Then

$$\frac{2}{N^2}E'_r(\theta) = \frac{1}{\beta}\frac{1}{g(x)^2}(f(x)g'(x) - f'(x)g(x))\Big(F(u) - \beta a\Big).$$

7

The next step is to understand the sign of the term $F(u) - \beta a$, which requires to understand the polynomials $P$ and $Q$. In the next section, we analyse the polynomial $P$ ($Q$ is explicit and much simpler). □

## 1.3 Concrete examples of small populations

**Example 1.1.** ($N = 3$). When $N = 3$, the total cost interference for the reward case is given by

$$E_r(\theta) = \frac{9\theta(1+e^x)}{2(1+e^x+e^{2x})}\left[(1+e^x)(1+1/2) + e^{2x}(e^{-x}+e^{-2x}/2)\right]$$
$$= \frac{9\theta}{4}\left(5 + \frac{4e^x - 1}{1 + e^x + e^{2x}}\right).$$

Note that since $x = \beta(\theta - a)$, $\frac{d}{d\theta}E(\theta) = \beta\frac{d}{dx}E(x)$.

$$\frac{4}{9}E'_r(\theta) = \frac{1}{\beta}\left(5 + \frac{4e^x - 1}{1 + e^x + e^{2x}}\right) + \theta\left(\frac{4e^x}{1+e^x+e^{2x}} - \frac{(4e^x - 1)(e^x + 2e^{2x})}{(1+e^x+e^{2x})^2}\right)$$
$$= \frac{1}{\beta}\left[\left(5 + \frac{4e^x - 1}{1+e^x+e^{2x}}\right) + (x + \beta a)\left(\frac{4e^x}{1+e^x+e^{2x}} - \frac{(4e^x-1)(e^x+2e^{2x})}{(1+e^x+e^{2x})^2}\right)\right]$$
$$= \frac{1}{\beta(1+u+u^2)^2}\left[(u+1)(u+4/5)(5u^2+5u+5) - (x+\beta a)u(4u^2-2u-5)\right] \quad (12)$$

where $a = c + \frac{b}{N-1} > 0$. Note that $x + \beta a = \beta\theta > 0$.

(1) If $x \leq \log\left(\frac{1+\sqrt{21}}{4}\right)$, then $4u^2 - 2u - 5 \leq 0$. It follows from (12) that $E'(x) > 0$ for all $\beta$. Hence $E$ is increasing for all $\beta$.

(2) If $x > \log\left(\frac{1+\sqrt{21}}{4}\right)$, then $4u^2 - 2u - 5 > 0$. Consider the function

$$F(u) := \frac{(u+1)(5u+4)(u^2+u+1)}{u(4u^2-2u-5)} - \log(u) \quad \text{for } u > \frac{1+\sqrt{21}}{4} =: u^0.$$

Then

$$EC'(x) = \frac{9}{4}\frac{u(4u^2-2u-5)}{\beta(1+u+u^2)^2}\Big[F(u) - (\beta a)\Big].$$

We have $\min_{u > u^0} F(u) = 10.9291 = f^*$ at $u^* = 4.29712$. Moreover, $f(u)$ is decreasing when $u < u^*$ and is increasing when $u > u^*$. The number of solutions of the equation $F(u) = \beta a$ and the sign of $E'(x)$ depends on the value of $\beta a$.

  (i) If $\beta a < f^*$, then $F(u) = \beta a$ has no solution and $E'(x) > 0$. Hence $E$ is increasing.

  (ii) If $\beta a = f^*$, then $F(u) = \beta a$ has a unique solution, and $E'(x) \geq 0$ and $E$ is non-decreasing.

  (iii) If $\beta a > f^*$, then $F(u) = \beta a$ has two solutions $(1+\sqrt{21})/4 < u_1 < u_2$.

$$E'_r(\theta) \begin{cases} > 0 & u < u_1, \\ < 0 & u_1 < u < u_2, \\ > 0 & u > u_2. \end{cases}$$

Thus $E_r$ is increasing when $u < u_1$, is decreasing when $u_1 < u < u_2$ and is increasing again when $u > u_2$.

**Algorithm for $N = 3$**

1. Compute

$$\delta = \begin{cases} -(c + \frac{b}{N-1}) & \text{PD games,} \\ -c\left(1 - \frac{r(N-n)}{n(N-1)}\right) & \text{PGG games} \end{cases} \quad (13)$$

2. Compute $\theta_0 = \frac{1}{(N-1)\beta} \log\left(\frac{\omega}{1-\omega}\right) - \delta$;

3. Compute $\beta^* = -\frac{f^*}{\delta} = -\frac{10.9291}{\delta}$.

4. If $\beta \leq \beta^*$ then
$$\theta = \theta_0, \quad \min E(\theta) = E(\theta_0).$$

5. If $\beta > \beta^*$
   (a) Compute $\theta_2 = \frac{\log u_2}{\beta} - \delta$ where $u_2$ is the largest zero of the equation $F(u) + \beta\delta = 0$.
   (b) If $\theta_0 \geq \theta_2$: $\theta^* = \theta_0, \quad \min E(\theta) = E(\theta_0)$;
   (c) if $\theta_0 < \theta_2$, compare $E(\theta_0)$ and $E(\theta_2)$: if $E(\theta_0) \leq E(\theta_2)$: $\theta^* = \theta_0, \quad \min E(\theta) = E(\theta_0)$; if $E(\theta_2) < E(\theta_0)$: $\theta^* = \theta_2, \quad \min E(\theta) = E(\theta_2)$.

**Example 1.2.** ($N = 4$) When $N = 4$, the total cost function is $E_r(\theta) = 8\theta\frac{f(x)}{g(x)}$, where

$$f(x) = (1 + e^x)\left[\frac{11}{6}(1 + e^x + e^{2x}) + e^{3x}(e^{-x} + e^{-2x}/2 + e^{-3x}/3)\right]$$
$$= (1 + e^x)(\frac{13}{6} + \frac{14}{6}e^x + \frac{17}{6}e^{2x})$$
$$= (1 + u)\left(13/6 + 14/6u + 17/6u^2\right),$$
$$g(x) = 1 + e^x + e^{2x} + e^{3x} = 1 + u + u^2 + u^3.$$

In this case,

$$F(u) = \frac{1}{2}\frac{(13 + 14u + 17u^2)(1 + u + u^2 + u^3)}{u(1+u)(7u^2 - 4u - 7)} - \log(u)$$

**Algorithm for $N = 4$**

1. Compute

$$\delta = \begin{cases} -(c + \frac{b}{N-1}) & \text{PD games,} \\ -c\left(1 - \frac{r(N-n)}{n(N-1)}\right) & \text{PGG games} \end{cases} \quad (14)$$

2. Compute $\theta_0 = \frac{1}{(N-1)\beta} \log\left(\frac{\omega}{1-\omega}\right) - \delta$;

3. Compute $u_0 > 0$ that solves $7u^2 - 4u - 7 = 0$; that is $u_0 = 3.4697$;

4. Compute $f^* = \min\{F(u) : u > u_0\}$, which gives $f^* = 6.64711$;

5. Compute $\beta^* = -\frac{f^*}{\delta}$.

6. If $\beta \leq \beta^*$ then
$$\theta = \theta_0, \quad \min E(\theta) = E(\theta_0).$$

7. If $\beta > \beta^*$
   (a) Compute $\theta_2 = \frac{\log u_2}{\beta} - \delta$ where $u_2$ is the largest zero of the equation $F(u) + \beta\delta = 0$.
   (b) if $E_r(\theta_0) \leq E_r(\theta_2)$: $\theta^* = \theta_0, \quad \min E_r(\theta) = E(\theta_0)$;
   (c) if $E_r(\theta_2) < E_r(\theta_0)$: $\theta^* = \theta_2, \quad \min E_r(\theta) = E_r(\theta_2)$.

## 1.4 The polynomial $P$

In this section, we will compute the coefficients of the polynomial $P$ and study its positive roots. To this end, we will need the following auxiliary lemma.

**Lemma 1.4.** *It holds that*

$$\sum_{j=m}^{n} \left( \frac{1}{N+j-k-2} + \frac{1}{N+j-k-1} - \frac{2}{N-1-j} \right) j$$
$$= 4(n+1-m) + (k+2-N) \sum_{j=N+m-k-2}^{N+n-k-2} \frac{1}{j} + (k-N+1) \sum_{j=N+m-k-1}^{N+n-k-1} \frac{1}{j} - 2(N-1) \sum_{j=N-1-n}^{N-1-m} \frac{1}{j}. \quad (15)$$

*In particular, we have*

$$\sum_{j=1}^{k} \left( \frac{1}{N+j-k-2} + \frac{1}{N+j-k-1} - \frac{2}{N-1-j} \right) j$$
$$= 4k + (k+2-N) \sum_{j=N-1-k}^{N-2} \frac{1}{j} + (k-N+1) \sum_{j=N-k}^{N-1} \frac{1}{j} - 2(N-1) \sum_{j=N-1-k}^{N-2} \frac{1}{j}. \quad (16)$$

*and*

$$\sum_{j=k-N+3}^{N-2} \left( \frac{1}{N+j-k-2} + \frac{1}{N+j-k-1} - \frac{2}{N-1-j} \right) j$$
$$= 4(2N-k-4) + (k+2-N) \sum_{j=1}^{2N-k-4} \frac{1}{j} + (k-N+1) \sum_{j=2}^{2N-k-3} \frac{1}{j} - 2(N-1) \sum_{j=1}^{2N-4-k} \frac{1}{j}$$
$$= 4(2N-k-4) - (4N-2k-5) \sum_{j=1}^{2N-k-4} \frac{1}{j} + (k-N+1) \left( \frac{1}{2N-k-3} - 1 \right). \quad (17)$$

*Proof.* By changing of variable $\bar{j} := N + j - k - 2$, we have

$$\sum_{j=m}^{n} \frac{j}{N+j-k-2} = \sum_{\bar{j}=N+m-k-2}^{N+n-k-2} \frac{k+2-N+\bar{j}}{\bar{j}}$$
$$= \sum_{\bar{j}=N+m-k-2}^{N+n-k-2} \left( 1 + \frac{k+2-N}{\bar{j}} \right)$$
$$= (n-m+1) + (k+2-N) \sum_{\bar{j}=N+m-k-2}^{N+n-k-2} \frac{1}{\bar{j}}.$$

Similarly we have

$$\sum_{j=m}^{n} \frac{j}{N+j-k-1} = \sum_{\bar{j}=N+m-k-1}^{N+n-k-1} \frac{k-N+1+\bar{j}}{\bar{j}}$$
$$= \sum_{\bar{j}=N+m-k-1}^{N+n-k-1} \left( 1 + \frac{k-N+1}{\bar{j}} \right)$$
$$= (n-m+1) + (k-N+1) \sum_{\bar{j}=N+m-k-1}^{N+n-k-1} \frac{1}{\bar{j}}.$$

and

$$\sum_{j=m}^{n} \frac{j}{N-1-j} = \sum_{\bar{j}=N-1-n}^{N-1-m} \frac{N-1-\bar{j}}{\bar{j}}$$

$$= \sum_{\bar{j}=N-1-n}^{N-1-m} \left(\frac{N-1}{\bar{j}} - 1\right)$$

$$= (N-1) \sum_{\bar{j}=N-1-n}^{N-1-m} \frac{1}{\bar{j}} - (n+1-m).$$

Thus

$$\sum_{j=m}^{n} \left(\frac{1}{N+j-k-2} + \frac{1}{N+j-k-1} - \frac{2}{N-1-j}\right) j$$

$$= 4(n+1-m) + (k+2-N) \sum_{j=N+m-k-2}^{N+n-k-2} \frac{1}{j}$$

$$+ (k-N+1) \sum_{j=N+m-k-1}^{N+n-k-1} \frac{1}{j} - 2(N-1) \sum_{j=N-1-n}^{N-1-m} \frac{1}{j}.$$

□

The next lemma provides a sharp estimates for the harmonic number that will be used later.

**Lemma 1.5.** *[Young, 1991] The harmonic number satisfies the followings estimates*

$$\log(n) + \gamma + \frac{1}{2n+1} \leq \sum_{j=1}^{n} \frac{1}{j} \leq \log(n) + \gamma + \frac{1}{2n-1}, \tag{18}$$

*where $\gamma = 0.5772...$ is the Euler–Mascheroni constant.*

In the next proposition we compute the coefficients of the polynomial $P$ in terms of the harmonic number. This proposition is the most technical result but is the key to the analysis of this paper.

**Proposition 1.6.** *Let $P(u)$ be the polynomial defined in (10). Then its coefficients are given by explicitly by*

$$p_k = \begin{cases} -(H_N + \frac{1}{N-2}) & k = 0, \\ 4(k+1) + (2k+4-4N) \sum_{j=N-1-k}^{N-1} \frac{1}{j} - \frac{k+1}{N-k-2} - H_N(k+1) & 1 \leq k \leq N-3, \\ 4(N-1) - 2NH_N & k = N-2, \\ 4N - 5 - \frac{1}{N-1} - NH_N & k = N-1, \\ 8N - 4k - 9 - (2N-k-1)H_N + (4N-2k-4) \sum_{j=2N-k-1}^{N-1} \frac{1}{j} - \frac{1}{2N-k-2} & N \leq k \leq 2N-5, \\ H_N + \frac{1}{2} & k = 2N-4. \end{cases} \tag{19}$$

*Furthermore, $P$ has exactly one positive root which is bigger than 1.*

*Proof.* The proof is rather technical and lengthy. To obtain the formula for $p_k$ we apply a general formula for coefficients of a product of two polynomials and use Lemma 1.4. Then we apply Decarte's rule of signs to show that $P$ has exactly one positive root.

11

Let $a_j := H_N + \frac{1}{N-1-j}$ for $j = 0, \ldots, N-2$. Suppose that $P(u) = \sum_{j=0}^{2N-4} p_j u^j$. Using the following formula of product of two polynomials:

$$\Big(\sum_{j=0}^{m} a_j x^j\Big)\Big(\sum_{j=0}^{n} b_j x^j\Big) = \sum_{k=0}^{m+n} \Big(\sum_{j=\max\{0,k-n\}}^{\min\{m,k\}} a_j b_{k-j}\Big) x^k,$$

we have

$(i)$ $\Big(\sum_{j=0}^{N-2}(H_N + \frac{1}{N-1-j})u^j\Big)\Big(\sum_{j=1}^{N-1} ju^{j-1}\Big) = \Big(\sum_{j=0}^{N-2} a_j u^j\Big)\Big(\sum_{j=1}^{N-2}(j+1)u^j\Big)$

$$= \sum_{k=0}^{2N-4}\Big(\sum_{j=\max(k-N+2,0)}^{\min(k,N-2)} a_j(k-j+1)\Big) u^k$$

$(ii)$ $\Big(\sum_{j=1}^{N-2}\Big(H_N + \frac{1}{N-1-j}\Big)ju^{j-1}\Big)\Big(\sum_{j=0}^{N-1} u^j\Big) = \Big(\sum_{j=0}^{N-3} a_{j+1}(j+1)u^j\Big)\Big(\sum_{j=0}^{N-1} u^j\Big)$

$$= \sum_{k=0}^{2N-4}\Big(\sum_{j=\max(0,k-N+1)}^{\min(k,N-3)} a_{j+1}(j+1)\Big) u^k.$$

$(iii)$ $\Big(\sum_{j=0}^{N-2}\Big(H_N + \frac{1}{N-1-j}\Big)u^j\Big)\Big(\sum_{j=0}^{N-1} u^j\Big) = \Big(\sum_{j=0}^{N-2} a_j u^j\Big)\Big(\sum_{j=0}^{N-1} u^j\Big)$

$$= \sum_{k=0}^{2N-3}\Big(\sum_{j=\max(0,k-N+1)}^{\min(N-2,k)} a_j\Big) u^k.$$

Hence

$(iii)$ $\Big[\Big(\sum_{j=0}^{N-2}(H_N + \frac{1}{N-1-j})u^j\Big)\Big(\sum_{j=1}^{N-1} ju^{j-1}\Big) - \Big(\sum_{j=1}^{N-2}\Big(H_N + \frac{1}{N-1-j}\Big)ju^{j-1}\Big)\Big(\sum_{j=0}^{N-1} u^j\Big)\Big](1+u)$

$$= \Big(\sum_{k=0}^{2N-4}\Big(\sum_{j=\max(k-N+2,0)}^{\min(k,N-2)} a_j(k-j+1) - \sum_{j=\max(0,k-N+1)}^{\min(k,N-3)} a_{j+1}(j+1)\Big) u^k\Big)(1+u)$$

$$= \Big(\sum_{k=0}^{2N-4} b_k u^k\Big)(1+u)$$

$$= \sum_{k=0}^{2N-3}\Big(\sum_{j=\max(0,k-1)}^{\min(k,2N-4)} b_j\Big) u^k$$

$$= b_0 + \sum_{k=1}^{2N-4}(b_{k-1} + b_k)u^k + b_{2N-4}u^{2N-3}$$

12

where for $k = 0, \ldots, 2N - 4$

$$b_k := \sum_{j=\max(k-N+2,0)}^{\min(k,N-2)} a_j(k-j+1) - \sum_{j=\max(0,k-N+1)}^{\min(k,N-3)} a_{j+1}(j+1) \qquad (20)$$

$$= \sum_{j=\max(k-N+3,1)}^{\min(k+1,N-1)} a_{k-j+1} j - \sum_{j=\max(1,k-N+2)}^{\min(k+1,N-2)} a_j j$$

$$= \begin{cases} \sum_{j=1}^{k+1}(a_{k-j+1} - a_j)j & k \le N-3, \\ \sum_{j=1}^{N-1} a_{k-j+1} j - \sum_{j=1}^{N-2} a_j j & k = N-2, \\ \sum_{j=k-N+3}^{N-1} a_{k-j+1} j - \sum_{j=k-N+2}^{N-2} a_j j & N-1 \le k \le 2N-4. \end{cases}$$

$$= \begin{cases} \sum_{j=1}^{k+1}(a_{k-j+1} - a_j)j & k \le N-3, \\ \sum_{j=1}^{N-2}(a_{k-j+1} - a_j)j + H_N(N-1) + 1 & k = N-2, \\ \sum_{j=k-N+3}^{N-2}(a_{k-j+1} - a_j)j + a_{k-N+2}(2N-3-k) & N-1 \le k \le 2N-5, \\ a_{N-2} & k = 2N-4. \end{cases}$$

Therefore, we obtain

$$P(u) = \sum_{k=0}^{2N-3} \Big( \sum_{j=\max(0,k-1)}^{\min(k,2N-4)} b_j - \sum_{j=\max(0,k-N+1)}^{\min(N-2,k)} a_j \Big) u^k,$$

that is for $k = 0, \ldots, 2N - 3$

$$p_k = \sum_{j=\max(0,k-1)}^{\min(k,2N-4)} b_j - \sum_{j=\max(0,k-N+1)}^{\min(N-2,k)} a_j \qquad (21)$$

$$= \begin{cases} b_0 - a_0 & k = 0, \\ b_{k-1} + b_k - \sum_{j=0}^{k} a_j & 1 \le k \le N-2, \\ b_{k-1} + b_k - \sum_{j=k-N+1}^{N-2} a_j & N-1 \le k \le 2N-4, \\ b_{2N-4} - a_{N-2} & k = 2N-3. \end{cases} \qquad (22)$$

In particular, the coefficient of $u^{2N-3}$ in $P(u)$ is

$$p_{2N-3} = b_{2N-4} - a_{N-2} = \Big(a_{N-2}(N-1) - a_{N-2}(N-2)\Big) - a_{N-2} = 0, \qquad (23)$$

thus $P(u)$ is a polynomial of degree $2N - 4$.

$$p_0 = b_0 - a_0 = (a_0 - a_1) - a_0 = -a_1 < 0.$$

**For $1 \le k \le N-3$:**

$$b_{k-1} + b_k = \sum_{j=1}^{k}(a_{k-j} - a_j)j + \sum_{j=1}^{k+1}(a_{k-j+1} - a_j)j$$

$$= \sum_{j=1}^{k}(a_{k-j} - 2a_j + a_{k-j+1})j + (a_0 - a_{k+1})(k+1)$$

$$= \sum_{j=1}^{k}\Big(\frac{1}{N+j-k-2} + \frac{1}{N+j-k-1} - \frac{2}{N-1-j}\Big)j + (\frac{1}{N-1} - \frac{1}{N-k-2})(k+1)$$

13

Therefore,

$$p_k = b_{k-1} + b_k - \sum_{j=0}^{k} a_j$$

$$= \sum_{j=1}^{k} \left( \frac{1}{N+j-k-2} + \frac{1}{N+j-k-1} - \frac{2}{N-1-j} \right) j + \left( \frac{1}{N-1} - \frac{1}{N-k-2} \right)(k+1)$$

$$- H_N(k+1) - \sum_{j=0}^{k} \frac{1}{N-1-j}$$

$$= \sum_{j=1}^{k} \left( \frac{1}{N+j-k-2} + \frac{1}{N+j-k-1} - \frac{2}{N-1-j} \right) j - H_N(k+1) - \sum_{j=0}^{k} \frac{1}{N-1-j}$$

$$- \frac{(k+1)^2}{(N-1)(N-k-2)}$$

$$= 4k + (k+2-N) \sum_{j=N-1-k}^{N-2} \frac{1}{j} + (k-N+1) \sum_{j=N-k}^{N-1} \frac{1}{j} - 2(N-1) \sum_{j=N-1-k}^{N-2} \frac{1}{j}$$

$$- H_N(k+1) - \sum_{j=0}^{k} \frac{1}{N-1-j} + \left( \frac{1}{N-1} - \frac{1}{N-k-2} \right)(k+1)$$

$$= 4(k+1) + (2k+4-4N) \sum_{j=N-1-k}^{N-1} \frac{1}{j} - \frac{k+1}{N-k-2} - H_N(k+1).$$

**For** $k = N-2$

$$b_{k-1} + b_k = b_{N-3} + b_{N-2}$$

$$= \sum_{j=1}^{N-2} (a_{k-j} - a_j) j + \sum_{j=1}^{N-2} (a_{k-j+1} - a_j) j + H_N(N-1) + 1$$

$$= \sum_{j=1}^{N-2} (a_{k-j} - 2a_j + a_{k-j+1}) j + H_N(N-1) + 1$$

$$= \sum_{j=1}^{N-2} \left( \frac{1}{N+j-k-2} + \frac{1}{N+j-k-1} - \frac{2}{N-1-j} \right) j + 1 + H_N(N-1)$$

$$= \sum_{j=1}^{N-2} \left( \frac{1}{j} + \frac{1}{1+j} - \frac{2}{N-1-j} \right) j + 1 + H_N(N-1)$$

$$= (N-2) + \sum_{j=1}^{N-2} \left( 1 - \frac{1}{1+j} \right) - 2 \sum_{j=1}^{N-2} \frac{j}{N-1-j} + H_N(N-1) + 1$$

$$= 2(N-2) - (H_N - 1) - 2[(N-1)H_N - (N-1)] + H_N(N-1) + 1$$

$$= 4(N-1) - N H_N,$$

where we have used the relation

$$\sum_{j=1}^{N-2} \frac{j}{N-1-j} = \sum_{j=1}^{N-2} \frac{N-1-j}{j} = (N-1) \sum_{j=1}^{N-2} \frac{1}{j} - (N-2) = (N-1)\left( H_N - \frac{1}{N-1} \right) - (N-2)$$

$$= (N-1)(H_N - 1). \tag{24}$$

14

We also have

$$\sum_{j=0}^{N-2} a_j = (N-1)H_N + \sum_{j=0}^{N-2} \frac{1}{N-1-j} = NH_N.$$

Thus

$$p_{N-2} = b_{N-3} + b_{N-2} - \sum_{j=0}^{N-2} a_j$$
$$= 4(N-1) - 2NH_N < 0,$$

where the last inequality follows from the fact that

$$H_N = \sum_{j=1}^{N-1} \frac{1}{j} > 1 + \frac{N-2}{N} = \frac{2(N-1)}{N}. \tag{25}$$

**For $k = N-1$:**

$$\begin{aligned}
b_{k-1} + b_k &= b_{N-2} + b_{N-1} \\
&= \sum_{j=1}^{N-2}(a_{N-1-j} - a_j)j + H_N(N-1) + 1 + \sum_{j=2}^{N-2}(a_{N-j} - a_j)j + a_1(N-2) \\
&= \sum_{j=2}^{N-2}\left(\frac{1}{j} + \frac{1}{j-1} - \frac{2}{N-1-j}\right)j + (1 - \frac{1}{N-2}) + H_N(2N-3) + 2 \\
&= \left(N - 3 + \sum_{j=1}^{N-3}\frac{j+1}{j}\right) - 2\sum_{j=2}^{N-2}\frac{j}{N-1-j} + H_N(2N-3) + 3 - \frac{1}{N-2} \\
&= 2(N-3) + \left(H_N - \frac{1}{N-2} - \frac{1}{N-1}\right) - 2\left((N-1)(H_N - 1) - \frac{1}{N-2}\right) + (2N-3)H_N + 3 - \frac{1}{N-2} \\
&= 4N - 5 - \frac{1}{N-1}
\end{aligned}$$

Thus

$$p_{N-1} = b_{N-2} + b_{N-1} - \sum_{k=0}^{N-2} a_j$$
$$= 4N - 5 - \frac{1}{N-1} - NH_N$$

15

**For $N \leq k \leq 2N - 5$:**

$$b_{k-1} + b_k = \sum_{j=k-N+2}^{N-2} (a_{k-j} - a_j)j + a_{k-N+1}(2N - 2 - k) + \sum_{j=k-N+3}^{N-2} (a_{k-j+1} - a_j)j + a_{k-N+2}(2N - 3 - k)$$

$$= \sum_{j=k-N+3}^{N-2} (a_{k-j+1} - 2a_j + a_{k-j})j + (a_{N-2} - a_{k-N+2})(k - N + 2)$$
$$+ a_{k-N+1}(2N - 2 - k) + a_{k-N+2}(2N - 3 - k)$$

$$= \sum_{j=k-N+3}^{N-2} \Big(\frac{1}{N+j-k-2} + \frac{1}{N+j-k-1} - \frac{2}{N-1-j}\Big)j + (H_N + 1)(k - N + 2)$$
$$+ (H_N + 1/(2N - k - 2))(2N - k - 2) + (3N - 2k - 5)(H_N + 1/(2N - k - 3))$$

$$= \sum_{j=k-N+3}^{N-2} \Big(\frac{1}{N+j-k-2} + \frac{1}{N+j-k-1} - \frac{2}{N-1-j}\Big)j$$
$$+ (4N - 2k - 5)H_N + k - N + 3 + \frac{3N - 2k - 5}{2N - k - 3}.$$

We also have

$$\sum_{j=k-N+1}^{N-2} a_j = (2N - k - 2)H_N + \sum_{j=k-N+1}^{N-2} \frac{1}{N-1-j}$$
$$= (2N - k - 2)H_N + \sum_{j=1}^{2N-k-2} \frac{1}{j}.$$

16

Thus using Lemma 1.4 we obtain

$$p_k = b_{k-1} + b_k - \sum_{j=k-N+1}^{N-2} a_j$$

$$= \sum_{j=k-N+3}^{N-2} \left( \frac{1}{N+j-k-2} + \frac{1}{N+j-k-1} - \frac{2}{N-1-j} \right) j + (4N-2k-5)H_N$$

$$+ k - N + 3 + \frac{3N-2k-5}{2N-k-3} - (2N-k-2)H_N - \sum_{j=1}^{2N-k-2} \frac{1}{j}$$

$$= \sum_{j=k-N+3}^{N-2} \left( \frac{1}{N+j-k-2} + \frac{1}{N+j-k} - \frac{2}{N-1-j} \right) j + (2N-k-3)H_N$$

$$+ k - N + 3 + \frac{3N-2k-5}{2N-k-3} - \sum_{j=1}^{2N-k-2} \frac{1}{j}$$

$$= 4(2N-k-4) - (4N-2k-5)\sum_{j=1}^{2N-k-4} \frac{1}{j} + (k-N+1)\left( \frac{1}{2N-k-3} - 1 \right)$$

$$+ (2N-k-3)H_N + k - N + 3 + \frac{3N-2k-5}{2N-k-3} - \sum_{j=1}^{2N-k-2} \frac{1}{j}$$

$$= 8N - 4k - 13 + (2N-k-3)H_N - (4N-2k-4)\sum_{j=1}^{2N-k-4} \frac{1}{j}$$

$$- \frac{2}{2N-k-3} - \frac{1}{2N-k-2}$$

$$= 8N - 4k - 9 - (2N-k-1)H_N + (4N-2k-4)\sum_{j=2N-k-1}^{N-1} \frac{1}{j} - \frac{1}{2N-k-2}$$

In particular, for $k = N$

$$p_N = 4N - 9 - (N-1)H_N + (2N-4)\frac{1}{N-1} - \frac{1}{N-2}.$$

Finally, for $k = 2N-4$

$$p_{2N-4} = b_{2N-5} + b_{2N-4} - (a_{N-3} + a_{N-2}) = 2a_{N-3} + a_{N-2} - (a_{N-3} + a_{N-2}) = a_{N-3} = H_N + \frac{1}{2}.$$

In conclusion, we obtain

$$p_k = \begin{cases} -(H_N + \frac{1}{N-2}) & k = 0, \\ 4(k+1) + (2k+4-4N)\sum_{j=N-1-k}^{N-1} \frac{1}{j} - \frac{k+1}{N-k-2} - H_N(k+1) & 1 \le k \le N-3, \\ 4(N-1) - 2NH_N & k = N-2, \\ 4N - 5 - \frac{1}{N-1} - NH_N & k = N-1, \\ 8N - 4k - 9 - (2N-k-1)H_N + (4N-2k-4)\sum_{j=2N-k-1}^{N-1} \frac{1}{j} - \frac{1}{2N-k-2} & N \le k \le 2N-5, \\ H_N + \frac{1}{2} & k = 2N-4. \end{cases}$$
(26)

Now we prove that there is only one change of signs in the sequence $\{p_k, 0 \le k \le 2N-4\}$.

Consider first $1 \le k \le N-3$. Since

$$\sum_{j=N-1-k}^{N-1} \frac{1}{j} \ge \frac{k+1}{N-1},$$

17

we have,

$$p_k \leq 4(k+1) - (4N - 2k - 4)\frac{k+1}{N-1} - \frac{k+1}{N-k-2} - H_N(k+1)$$
$$\leq (k+1)\left(4 - \frac{4N - 2k - 4}{N-1} - \frac{1}{N-k-2} - H_N\right)$$
$$= (k+1)\left(\frac{2k}{N-1} - H_N - \frac{1}{N-k-2}\right) < 0,$$

where to obtain the last inequality we have used (25) so that

$$H_N \geq \frac{2(N-1)}{N} > \frac{2(N-3)}{N-1} \geq \frac{2k}{N-1}.$$

Thus $p_k < 0$ for all $0 \leq k \leq N - 3$. Next we consider $p_m =$ with $N \leq m \leq 2N - 5$. It will be more convenient to set $m = 2N - 4 - k$ and consider instead $p_{2N-4-k}$. We have

$$p_{2N-4-k} = 8N - 4(2N - 4 - k) - 9 - (2N - (2N - 4 - k) - 1)H_N$$
$$+ (4N - 2(2N - 4 - k) - 4)\sum_{j=k+3}^{N-1}\frac{1}{j} - \frac{1}{k+2}.$$
$$= 4k + 7 - (k+3)H_N + (2k+4)\sum_{j=k+3}^{N-1}\frac{1}{j} - \frac{1}{k+2}$$
$$= 4k + 7 - (k+3)H_N + (2k+4)(H_N - H_{k+3}) - \frac{1}{k+2}$$
$$= 4k + 7 + (k+1)H_N - (2k+4)H_{k+3} - \frac{1}{k+2}.$$

By lemma 1.5

$$\log(k+2) + \gamma + \frac{1}{2(k+2)+1} \leq H_{k+3} = \sum_{j=1}^{k+2}\frac{1}{j} \leq \log(k+2) + \gamma + \frac{1}{2(k+2)-1},$$

we have

$$4k + 7 + (k+1)H_N - (2k+4)(\log(k+2) + \gamma + \frac{1}{2k+3}) - \frac{1}{k+2} \leq p_{2N-k-4}$$
$$\leq 4k + 7 + (k+1)H_N - (2k+4)(\log(k+2) + \gamma + \frac{1}{2k+5}) - \frac{1}{k+2} \quad (27)$$

Consider the equation

$$4x + 7 + (x+1)H_N - 2(x+2)(\log(x+2) + \gamma + \frac{1}{2x+5}) - \frac{1}{x+2} = 0 \text{ for } x \geq 0, \quad (28)$$

which is equivalent to

$$h(x) := 2\log(x+2) + \frac{1}{x+2} + \frac{H_N}{x+2} + \frac{1}{(x+2)^2} + \frac{2}{2x+5} = 4 + H_N - 2\gamma.$$

We have

$$h'(x) = \frac{8x^4 + (64 - 4H_N)x^3 + (182 - 28H_N)x^2 + (207 - 65H_N)x + (68 - 50H_N)}{(x+2)^3(2x+5)^2}.$$

Let $r(x) := 8x^4 + (64 - 4H_N)x^3 + (182 - 28H_N)x^2 + (207 - 65H_N)x + (68 - 50H_N) := r_4x^4 + r_3x^3 + r_2x^2 + r_1x + r_0$, where $r_4 := 8, r_3 = 64 - 4H_N, r_2 = 182 - 28H_N, r_1 = 207 - 65H_N, r_0 := 68 - 50H_N$, be the enumerator of $h'(x)$. Since $N \geq 2$, $H_N \geq 1.5$, it follows that $r_0 < 0$. There is only one change of signs in the coefficients of $r(x)$ since:

(i) if $H_N \leq 207/65$ then $r_0 < 0 < r_1, r_2, r_3, r_4$,

(ii) if $207/65 < H_N < 182/28$, then $r_0, r_1 < 0$ while $r_2, r_3, r_4 > 0$,

(iii) if $182/28 < H_N \leq 16$, then $r_0, r_1, r_2 < 0$ while $r_3, r_4 > 0$,

(iv) $16 < H_N$ then $r_0, r_1, r_2, r_3 < 0$, while $r_4 > 0$.

Thus, by Decarte's rule of signs, $r(x)$ has a unique positive solution $x_0 > 0$. It follows that $h'(x) < 0$ when $x < x_0$ and $h'(x) > 0$ when $x > x_0$. Hence $h(x)$ is decreasing when $x < x_0$, increasing when $x > x_0$ and attains a global minimum at $x = x_0$. Since $f(0) = 2\log(2) + 1/2 + H_N/2 + 1/4 + 2/5 < 4 + H_N - 2\gamma$, Equation (28) has a unique solution $x^* > x_0$. If $k \geq x^* > x_0$ then, since $h(x)$ is increasing for $x > x_0$, $4 + H_N - 2\gamma \leq h(k) = 2\log(k+2) + \frac{1}{k+2} + \frac{H_N}{k+2} + \frac{1}{(k+2)^2} + \frac{2}{2k+5}$, that is

$$4k + 7 + (k+1)H_N - (2k+4)(\log(k+2) + \gamma + \frac{1}{2k+5}) - \frac{1}{k+2} \leq 0.$$

Hence if $k > x^*$ then

$$p_{2N-4-k} < 4k + 7 + (k+1)H_N - (2k+4)(\log(k+2) + \gamma + \frac{1}{2k+5}) - \frac{1}{k+2} < 0.$$

Similarly, let

$$\bar{g}(x) := 2\log(x+2) + \frac{1}{x+2} + \frac{H_N}{x+2} + \frac{1}{(x+2)^2} + \frac{2}{2x+3},$$

and consider the equation

$$\bar{g}(x) = 4 + H_N - 2\gamma. \tag{29}$$

$\bar{g}'(x) = 0$ has a unique positive solution; $\bar{g}(x)$ is decreasing when $x < x_1$, increasing when $x > x_1$ and attains a global minimum at $x = x_1$; equation (29) has a unique solution $\hat{x}$.

If $k \leq \hat{x}$ then $h(k) \leq 4 + H_N - 2\gamma$, that is

$$4k + 7 + (k+1)H_N - (2k+4)(\log(k+2) + \gamma + \frac{1}{2k+5}) - \frac{1}{k+2} \geq 0.$$

Hence if $k \leq \hat{x}$ then $p_{2N-k-4} > 0$. Since

$$\bar{g}(\hat{x}) = h(x^*) < \bar{g}(x^*)$$

and $\bar{g}$ is increasing in $(x_1, \infty)$, we have $\hat{x} < x^*$. We now prove that they are closed together.

$$\bar{g}(x^*) - h(x^*) = \bar{g}(x^*) - \bar{g}(\hat{x}) + \bar{g}(\hat{x}) - h(x^*) = \bar{g}(x^*) - \bar{g}(\hat{x}).$$

Hence

$$\frac{4}{(2x^*+3)(2x^*+5)} = 2(\log(x^*+2) - \log(\hat{x}+2)) + (H_N+1)\Big(\frac{1}{x^*+2} - \frac{1}{\hat{x}+2}\Big)$$
$$+ \Big(\frac{1}{(x^*+2)^2} - \frac{1}{(\hat{x}+2)^2}\Big) + \Big(\frac{2}{2x^*+3} - \frac{2}{2\hat{x}+3}\Big)$$
$$= (x^* - \hat{x})\Big(\frac{2}{\xi+2} - \frac{H_N+1}{(x^*+2)(\hat{x}+2)} - \frac{1}{(\hat{x}+2)^2(x^*+2)} - \frac{1}{(\hat{x}+2)(x^*+2)^2}$$
$$- \frac{4}{(2x^*+3)(2\hat{x}+3)}\Big), \tag{30}$$

where $\hat{x} < \xi < x^*$. Since

$$\bar{g}(H_N - 1) = 2\log(H_N+1) + 1 + \frac{1}{H_N+1} + \frac{2}{2H_N+1} < 4 + H_N - 2\gamma = \bar{g}(\hat{x}),$$

19

we have $m := H_N - 1 < \hat{x} < x^*$. For $N \geq 5$, we have $x^* > \hat{x} > m \geq \frac{77}{60}$. Hence

$$\frac{H_N+1}{(\hat{x}+2)(x^*+2)} + \frac{1}{(\hat{x}+2)^2(x^*+2)} + \frac{1}{(\hat{x}+2)(x^*+2)^2} + \frac{4}{(2x^*+3)(2\hat{x}+3)}$$
$$< \frac{1}{x^*+2} + \frac{1}{(m+2)^2(x^*+2)} + \frac{1}{(m+2)^2(x^*+2)} + \frac{4}{(2m+3)(2x^*+3)}.$$

It follows that

$$\frac{2}{\xi+2} - \frac{H_N+1}{(x^*+2)(\hat{x}+2)} - \frac{1}{(\hat{x}+2)^2(x^*+2)} - \frac{1}{(\hat{x}+2)(x^*+2)^2} - \frac{4}{(2x^*+3)(2\hat{x}+3)}$$
$$> \frac{2}{x^*+2} - \Big(\frac{1}{x^*+2} + \frac{1}{(m+2)^2(x^*+2)} + \frac{1}{(m+2)^2(x^*+2)} + \frac{4}{(2m+3)(2x^*+3)}\Big)$$
$$= \frac{(4m^3+18m^2+16m-4)x^* + (6m^3+25m^2+16m-14)}{(m+2)^2(2m+3)(x^*+2)(2x^*+3)}. \qquad (31)$$

From (30) and (31) we deduce

$$0 < (x^* - \hat{x}) < \frac{4(m+2)^2(2m+3)(x^*+2)}{(2x^*+5)(4m^3+18m^2+16m-4)x^* + (6m^3+25m^2+16m-14)}$$
$$< \frac{2(m+2)^2(2m+3)}{(4m^3+18m^2+16m-4)x^* + (6m^3+25m^2+16m-14)}$$
$$< \frac{2(m+2)^2(2m+3)}{(4m^3+18m^2+16m-4)m + (6m^3+25m^2+16m-14)}$$
$$< 1,$$

where the last inequality is because $m > 1.28$.

Let $k_0$ be the position of change of signs in the coefficients $\{p_k, 0 \leq k \leq 2N-4\}$, that is $p_k < 0$ for $0 \leq k \leq k_0$ and $p_k > 0$ for $k_0 + 1 \leq k \leq 2N - 4$. Since $p_N > 0$ for $N \geq 27$, we have

$$k_0 = \begin{cases} N-2 & N \leq 25, \\ N-1 & N \in \{26, 27\}, \\ 2N-4-\lfloor 2N-x^*-4 \rfloor & N \geq 27. \end{cases} \qquad (32)$$

Since there is only one change of signs in the sequence of coefficients of $P$, by Decarte's rule of signs, $P$ has a unique positive root. This root is bigger than 1 since we have

$$P(1) = 2\Big[\Big(\sum_{j=0}^{N-2}(H_N + \frac{1}{N-1-j})\Big)\Big(\sum_{j=1}^{N-1} j\Big) - N\sum_{j=1}^{N-2}\Big(H_N + \frac{1}{N-1-j}\Big)j\Big]$$
$$- N\sum_{j=0}^{N-2}\Big(H_N + \frac{1}{N-1-j}\Big)$$
$$= \Big[(N-1)N\Big((N-1)H_N + \sum_{j=0}^{N-2}\frac{1}{N-1-j}\Big) - N\Big((N-2)(N-1)H_N + 2\sum_{j=1}^{N-2}\frac{N-1-j}{j}\Big)\Big]$$
$$- N\Big((N-1)H_N + \sum_{j=0}^{N-2}\frac{1}{N-1-j}\Big)$$
$$= \Big[(N-1)N^2 H_N - N((N-2)(N-1)H_N + 2(N-1)H_N - 2 - 2(N-2))\Big] - N^2 H_N$$
$$= 2N(N-1) - N^2 H_N < 0$$

where the first inequality follows from (25). □

## 1.5 The function $F$

Next we study the function $F$ defined in (11). We will need the following lemma, which establishes a relation between the polynomials $P$ and $Q$.

**Lemma 1.7.** *For $u > 0$, we have*

$$Q(u) \geq (1+u)uP(u). \tag{33}$$

*As a consequence, for $\beta \leq \frac{1}{a}$ we have $\frac{d}{d\theta} E_r(\theta) \geq 0$.*

*Proof.* The statement (33) is equivalent to

$$\sum_{j=0}^{N-2} \left( H_N + \frac{1}{N-1-j} \right) u^j \Big) \Big( \sum_{j=0}^{N-1} u^j \Big) = \sum_{k=0}^{2N-3} \Big( \sum_{j=\max(0,k-N+1)}^{\min(N-2,k)} a_j \Big) u^k \geq uP(u), \tag{34}$$

which we now prove. Let

$$h_k := \sum_{j=\max(0,k-N+1)}^{\min(N-2,k)} a_j = \begin{cases} \sum_{j=0}^{k} a_j & 0 \leq k \leq N-2, \\ \sum_{j=k-N+1}^{N-2} a_j & N-1 \leq k \leq 2N-3 \end{cases}$$

$$= \begin{cases} (k+1)H_N + \sum_{j=N-1-k}^{N-2} \frac{1}{j} & 0 \leq k \leq N-2, \\ (2N-k-2)H_N + \sum_{j=1}^{2N-k-2} \frac{1}{j} & N-1 \leq k \leq 2N-3 \end{cases}$$

Now we show that $h_{k+1} \geq p_k$ for all $k = 0, \ldots, 2N-4$.

For $k = 0$: $h_1 > 0 > p_0$.

For $1 \leq k \leq N-3$, we have

$$h_{k+1} = (k+2)H_N + \sum_{j=N-k-2}^{N-2} \frac{1}{j}, \quad p_k = 4(k+1) + (2k+4-4N) \sum_{j=N-1-k}^{N-1} \frac{1}{j} - \frac{k+1}{N-k-2} - H_N(k+1).$$

Since $H_N \geq 2$ ($N \geq 4$), $(2k+3)H_N \geq 4(k+1)$. It follows that $h_{k+1} \geq p_k$.

For $k = N-2$

$$h_{N-1} = (N-1)H_N + \sum_{j=1}^{N-1} \frac{1}{j} = NH_N > 4(N-1) - 2NH_N$$

since according to (25), $3NH_N > 6(N-1) > 4(N-1)$.

For $k = N-1$,

$$h_N = (N-2)H_N + \sum_{j=1}^{N-2} \frac{1}{j} = (N-1)H_N - \frac{1}{N-1} > 4N - 5 - \frac{1}{N-1} - NH_N = p_{N-1}$$

since $H_N \geq 2$ ($N \geq 4$), $(2N-1)H_N \geq 4N - 2 > 4N - 5$.

For $N \leq k \leq 2N-5$. We have

$$h_{k+1} = (2N-k-3)H_N + \sum_{j=1}^{2N-k-3} \frac{1}{j},$$

$$p_k = 8N - 4k - 9 - (2N-k-1)H_N + (4N-2k-4) \sum_{j=2N-k-1}^{N-1} \frac{1}{j} - \frac{1}{2N-k-2}$$

$$= 8N - 4k - 9 - (2N-k-1)H_N + (4N-2k-4)[H_N - H_{2N-k-1}] - \frac{1}{2N-k-2}.$$

21

Thus, we obtain
$$h_{k+1} - p_k = (4N - 2k - 3)H_{2N-k-1} - (8N - 4k - 9).$$
Since $2N - k - 1 \geq 4$, we have $H_{2N-k-1} \geq 2$, thus
$$h_{k+1} - p_k \geq 2(4N - 2k - 3) - (8N - 4k - 9) = 3 > 0,$$
i.e., $h_{k+1} > p_k$ for all $N \leq k \leq 2N - 5$

Finally for $k = 2N - 4$, we have
$$h_{2N-3} = H_N + 1 > H_N + \frac{1}{2} = p_{2N-4}.$$
In conclusion, we have proved that $h_{k+1} > p_k$ for all $k = 0, \ldots, 2N - 4$.

As a consequence, if $\beta a \leq 2$, then
$$0 < x + \beta a \leq x + 2 \leq e^x + 1 = 1 + u,$$
thus
$$(x + \beta a)(f(x)g'(x) - f'(x)g(x)) = (x + \beta a)uP(u) \leq (1 + u)uP(u) \leq Q(u) = f(x)g(x),$$
which implies that $E'(x) \geq 0$. This completes the proof of the lemma. $\square$

Let $u^0 > 1$ be the unique positive root of $P$ obtained in Proposition 1.6. Then
$$\{u \in \mathbb{R} : P(u) > 0\} = \{u > u_0\}.$$
Lemma 1.7 provides an upper bound $\bar{\beta}$ of $\beta$ such that $E(\theta) \geq 0$ for all $\beta \leq \bar{\beta}$. To obtain the optimal bound we need to study the minimization problem
$$\min\{F(u) \mid u > u^0\}.$$

**Proposition 1.8.** *$F$ has a global minimizer on $(u^0, +\infty)$. Furthermore, for $N \leq N_0 = 100$, $F$ has a unique minimizer $u^*$ and $F$ is decreasing when $u^0 < u < u^*$ and is increasing when $u > u^*$.*

*Proof.* It follows from Lemma 1.7 that $F$ is bounded below on $(u^0, +\infty)$. Further more, it is continuous and coercive in this interval since
$$\lim_{u \to (u^0)^+} F(u) = \lim_{u \to +\infty} F(u) = +\infty.$$
Therefore $F$ has a global minimizer on $(u^0, +\infty)$. We have
$$F'(u) = \frac{uQ'(u)P(u) - Q(u)(P(u) + uP'(u))}{u^2 P(u)^2} - \frac{1}{u} =: \frac{M(u)}{u^2 P(u)^2},$$
where
$$M(u) := uQ'(u)P(u) - Q(u)(P(u) + uP'(u)) - uP(u)^2. \tag{35}$$
Then $M$ is a polynomial of degree $4N - 6$. The leading coefficient $m_{4N-6}$ of $M$ is
$$m_{4N-6} = q_{2N-2}p_{2N-4} > 0.$$
The sign of $F'$ is the same as that of $M$. Thus the number of changes of signs of $F'$ is equal to the number of positive roots of $M$. Since $Q(0) > 0$ while $P(0) < 0$, we have
$$M(0) = -Q(0)P(0) > 0$$
In addition, since $P$ is a polynomial whose coefficient of the highest degree is positive and $u^0$ is the unique positive root of $P$, it follows that $P'(u^0) > 0$. Hence
$$M(u^0) = -u^0 Q(u^0) P'(u^0) < 0.$$

Since $\lim_{u \to +\infty} M(u) = +\infty$, $M$ has at least two positive roots, one is belong to $(0, u^0)$ and one is in $(u^0, +\infty)$. To determine the number of positive roots of $M$ we can use Sturm's theorem. To this end, we construct the Sturm sequence of $M$ as follows

$$M_0 = M,$$
$$M_1 = M',$$
$$M_{i+1} = -\text{rem}(M_{i-1}, M_i),$$

for $i \geq 1$, where $M'$ is the derivative of $M$ and $\text{rem}(M_{i-1}, M_i)$ is the remainder of the Euclidean division of $M_{i-1}$ by $M_i$. The number of changes of signs at $\xi$ of the Sturm sequence of $M$, denoted by $s(\xi)$, is the number of sign changes, by ignoring zeros, in the sequence of real numbers

$$M_0(\xi), M_1(\xi), \ldots$$

According to Sturm theorem, the number of distinct positive real roots of $M$ is $s(0) - s(+\infty)$, where the sign at $+\infty$ of a polynomial is the sign of its leading coefficient. Sturm' algorithm is rather analytically intricate but is easily implemented using mathematical softwares such as Mathematica.

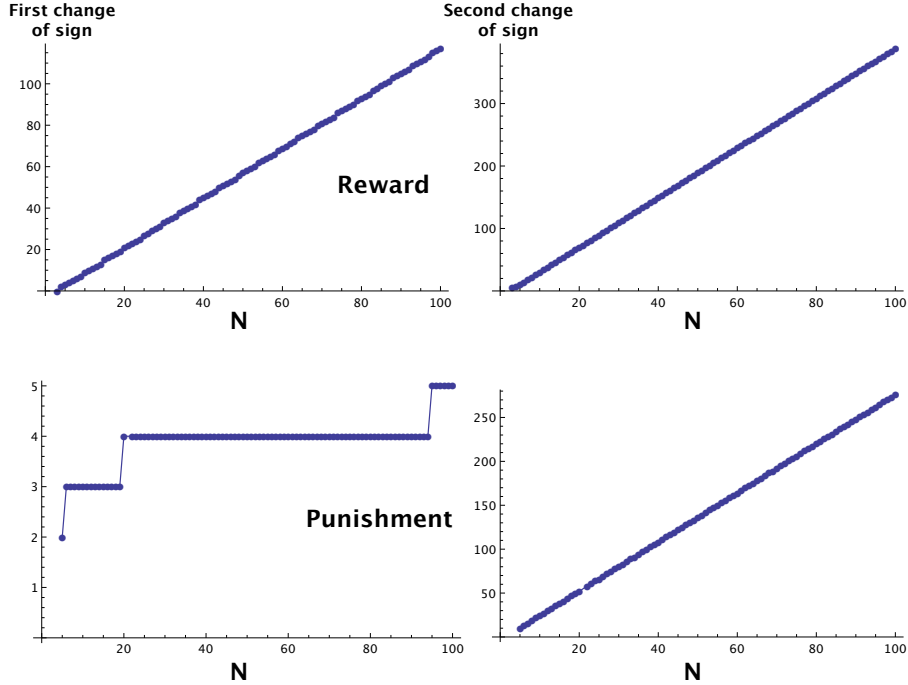

Figure 4: **Numerical results showing there are always two changes of sign in the sequence of coefficients of M, for $N$ up to 100.** We show the positions at which the changes of signs occur, for reward (top row) and punishment (bottom row).

We conjecture that $M$ actually has exactly two roots. One possible way to prove this conjecture is to show that there are two changes of signs in the sequence of coefficients of $M$. This is because if there are two changes of signs in the sequence of coefficients of $M$ then $M$ has either 2 or none positive roots according to Descartes' rule of signs, which follows that $M$ has exactly 2 positive roots since we already proved above that $M$ has at least two positive roots. However, finding the signs of coefficients of $M$ is analytically challenging since formulas for the coefficients of $M$ is extremely intricate. By direct computations (using Mathematica) we verify that this conjecture is true for $N \leq N_0 = 100$. As a consequence, for $N \leq N_0$, $F'$ has only one zero, denoted by $u^*$, on $(u^0, +\infty)$. Then $u^*$ is the unique minimizer of $F$ on $(u^0, +\infty)$. Furthermore, $F$ is decreasing on $(u^0, u^*)$ and is increasing on $(u^*, +\infty)$. □

## 1.6 Finite population estimates

In this section, we establish the finite population estimates for $E_r(\theta)$. To this end, we need the following auxiliary lemma.

**Lemma 1.9.** *For all $x \in \mathbb{R}$, we have*

$$0 \leq \frac{e^x + \ldots + e^{(N-2)x}}{1 + e^x + \ldots + e^{(N-1)x}} \leq \frac{N-2}{N}. \tag{36}$$

*Proof.* Let $u = e^x$ and $z(u) := \frac{u + \ldots + u^{N-2}}{1 + u + \ldots + u^{N-1}}$. Then

$$z(u) = \begin{cases} \frac{N-2}{N} & u = 1, \\ \frac{u^{N-1} - u}{u^N - 1} & u \neq 1. \end{cases}$$

For $u \neq 1$, we have

$$z'(u) = \frac{(N-1)u^N - u^{2N-2} - (N-1)u^{N-2} + 1}{(u^N - 1)^2}.$$

To find the roots of $f'(u)$ we need to factorize the enumerator in the above expression.

$$(N-1)u^N - u^{2N-2} - (N-1)u^{N-2} + 1$$
$$= (1 - u^{2N-2}) - (N-1)u^{N-2}(1-u)$$
$$= (1-u)(1 + u + \ldots + u^{2N-3} - (N-1)u^{N-2} - (N-1)u^{N-1})$$
$$= (1-u)\Big[(1 - u^{N-2}) + u(1 - u^{N-3}) + \ldots + u^{N-3}(1-u) - u^{N-1}(1-u)$$
$$\qquad - u^{N-1}(1 - u^2) - \ldots - u^{N-1}(1 - u^{N-2})\Big]$$
$$= (1-u)\Big[(1-u)(u^{N-3} - u^{N-1}) + (1-u^2)(u^{N-4} - u^{N-1}) + \ldots + (1-u^{N-2})(1 - u^{N-1})\Big]$$
$$= (1-u)\Big[(1-u)u^{N-3}(1-u^2) + (1-u^2)u^{N-4}(1-u^3) + \ldots + (1-u^{N-2})(1-u^{N-1})\Big]$$
$$= (1-u)^3\Big[u^{N-3}(1+u) + u^{N-4}(1+u)(1+u+u^2) + \ldots + (1+u+\ldots+u^{N-3})(1+u+\ldots+u^{N-2})\Big].$$

It follows that $z'(u) > 0$ in $(0,1)$ and $z'(u) < 0$ in $(1, +\infty)$. Hence $\max z(u) = f(1) = \frac{N-2}{N}$. In addition, we have

$$z(0) = 0 = \lim_{u \to +\infty} f(u).$$

Thus $0 \leq z(u) \leq \frac{N-2}{N}$ for all $u > 0$, which is the required statement. $\square$

**Lemma 1.10.** *It holds that*

$$\frac{N^2 \theta}{2}\left(H_N + \frac{1}{N-1}\right) \leq E_r(\theta) \leq N(N-1)\theta\Big(H_N + 1\Big). \tag{37}$$

*Proof.* We have

$$\frac{(1+e^x)(1 + e^x + \ldots + e^{(N-2)x})}{1 + e^x + \ldots + e^{(N-1)x}} = 1 + \frac{e^x + \ldots + e^{(N-2)x}}{1 + e^x + \ldots + e^{(N-1)x}}.$$

Using Lemma 1.9, we get

$$1 \leq \frac{(1+e^x)(1 + e^x + \ldots + e^{(N-2)x})}{1 + e^x + \ldots + e^{(N-1)x}} \leq \frac{2(N-1)}{N}.$$

Since

$$\frac{1}{N-1}\sum_{j=0}^{N-2} e^{jx} \leq e^{(N-1)x}\sum_{j=1}^{N-1} \frac{e^{-jx}}{j} = \sum_{j=0}^{N-2} \frac{e^{jx}}{N-1-j} \leq \sum_{j=0}^{N-2} e^{jx}.$$

24

we obtain the following estimates

$$\frac{1}{N-1} \le \frac{1}{N-1}\frac{(1+e^x)(1+e^x+\ldots+e^{(N-2)x})}{1+e^x+\ldots+e^{(N-1)x}} \le \frac{(1+e^x)}{1+e^x+\ldots+e^{(N-1)x}}e^{(N-1)x}\sum_{j=1}^{N-1}\frac{e^{-jx}}{j}$$

$$\le \frac{(1+e^x)(1+e^x+\ldots+e^{(N-2)x})}{1+e^x+\ldots+e^{(N-1)x}} \le \frac{2(N-1)}{N}.$$

Thus for $\theta > 0$ we have

$$\frac{N^2\theta}{2}\Big(\sum_{j=1}^{N-1}\frac{1}{j} + \frac{1}{N-1}\Big) \le E_r(\theta) \le \frac{N^2\theta}{2}\frac{2(N-1)}{N}\Big(\sum_{j=1}^{N-1}\frac{1}{j} + 1\Big) = N(N-1)\theta\Big(\sum_{j=1}^{N-1}\frac{1}{j} + 1\Big).$$

As a consequence, since $E$ is smooth as a function of $\theta$, by the extreme value theorem, it attains a minimum in any bounded interval $[\theta_0, \theta_1]$ where $\theta_0 \ge 0$. □

## 1.7 Asymptotic limits

In this section, we will establish the asymptotic limits of $E_r(\theta)$ as $N \to +\infty$, $\beta \to 0$ and $\beta \to +\infty$ where $N$ is the population size and $\beta$ is the selection intensity.

**Proposition 1.11** (Infinite population limit). *We have*

$$\lim_{N\to+\infty} \frac{E_r(\theta)}{\frac{N^2\theta}{2}(\ln N + \gamma)} = \begin{cases} 1 + e^{-\beta|\theta-c|} & \text{for DG,} \\ 1 + e^{-\beta|\theta-c|}e^{\beta c\frac{r}{n}} & \text{for PGG.} \end{cases} \quad (38)$$

*Proof.*

$$\frac{N^2\theta}{2}\frac{(1+e^x)(1+\ldots+e^{(N-2)x})}{1+\ldots+e^{(N-1)x}}\Big(\sum_{j=1}^{N-1}\frac{1}{j} + \frac{1}{N-1}\Big) \le E_r(\theta)$$

$$\le \frac{N^2\theta}{2}\frac{(1+e^x)(1+\ldots+e^{(N-2)x})}{1+\ldots+e^{(N-1)x}}\Big(\sum_{j=1}^{N-1}\frac{1}{j} + 1\Big)$$

Hence

$$\frac{(1+e^x)(1+\ldots+e^{(N-2)x})}{1+\ldots+e^{(N-1)x}}\frac{\big(\sum_{j=1}^{N-1}\frac{1}{j} + \frac{1}{N-1}\big)}{\ln N + \gamma} \le \frac{E_r(\theta)}{\frac{N^2\theta}{2}(\ln N + \gamma)}$$

$$\le \frac{(1+e^x)(1+\ldots+e^{(N-2)x})}{1+\ldots+e^{(N-1)x}}\frac{\big(\sum_{j=1}^{N-1}\frac{1}{j} + 1\big)}{\ln N + \gamma} \quad (39)$$

We recall that $x = \beta(\theta + \delta)$, where

$$\delta = \delta(N) = \begin{cases} -(c + \frac{b}{N-1}) & \text{for DG,} \\ -c(1 - \frac{r(N-n)}{n(N-1)}) & \text{for PGG.} \end{cases}$$

Let

$$w(\theta) := \frac{(1+e^x)(1+\ldots+e^{(N-2)x})}{1+\ldots+e^{(N-1)x}}.$$

Then we have

$$w(\theta) = 1 + \frac{e^x + \ldots e^{(N-2)x}}{1+\ldots+e^{(N-1)x}} = \begin{cases} \frac{2(N-1)}{N} & \text{if } x = 0, \\ 1 + e^x\frac{1-e^{(N-3)x}}{1-e^{(N-1)x}} & \text{if } x \ne 0. \end{cases}$$

Since for fixed $\theta$, $x$ equals to 0 for only one value of $N$, we can consider $x \neq 0$ when we study the limit $N \to +\infty$. We have

$$\lim_{N\to+\infty} e^x = \lim_{N\to+\infty} e^{\beta(\theta+\delta(N))} = \begin{cases} e^{\beta(\theta-c)} & \text{for DG,} \\ e^{\beta(\theta-c(1-\frac{r}{n}))} & \text{for PGG.} \end{cases}$$

For DG:

$$\lim_{N\to+\infty} \frac{1-e^{(N-3)x}}{1-e^{(N-1)x}} = \lim_{N\to+\infty} \frac{1-e^{-\beta b(N-3)/(N-1)}e^{\beta(\theta-c)(N-3)}}{1-e^{-\beta b}e^{\beta(\theta-c)(N-1)}} = \begin{cases} 1 & \text{if } \theta \leq c, \\ e^{-2\beta(\theta-c)} & \text{if } \theta > c. \end{cases}$$

For PGG:

$$\lim_{N\to+\infty} \frac{1-e^{(N-3)x}}{1-e^{(N-1)x}} = \lim_{N\to+\infty} \frac{1-e^{(N-3)\beta cr\frac{(N-n)}{n(N-1)}}e^{\beta(\theta-c)(N-3)}}{1-e^{\beta cr\frac{(N-n)}{n}}e^{\beta(\theta-c)(N-1)}} = \begin{cases} 1 & \text{if } \theta \leq c, \\ e^{-2(\theta-c)} & \text{if } \theta > c. \end{cases}$$

It follows that for DG

$$\lim_{N\to+\infty} w(\theta) = \begin{cases} 1+e^{\beta(\theta-c)} & \text{if } \theta \leq c, \\ 1+e^{-\beta(\theta-c)} & \text{if } \theta > c \end{cases}$$
$$= 1+e^{-\beta|\theta-c|}. \tag{40}$$

and for PGG

$$\lim_{N\to+\infty} w(\theta) = \begin{cases} 1+e^{\beta(\theta-c)}e^{\beta c\frac{r}{n}} & \text{if } \theta \leq c, \\ 1+e^{-\beta(\theta-c)}e^{\beta c\frac{r}{n}} & \text{if } \theta > c \end{cases}$$
$$= 1+e^{-\beta|\theta-c|}e^{\beta c\frac{r}{n}}. \tag{41}$$

From Lemma 1.5, we deduce that

$$\lim_{N\to+\infty} \frac{\left(\sum_{j=1}^{N-1}\frac{1}{j}+\frac{1}{N-1}\right)}{\ln N+\gamma} = \lim_{N\to+\infty} \frac{\left(\sum_{j=1}^{N-1}\frac{1}{j}+1\right)}{\ln N+\gamma} = 1 \tag{42}$$

From (39), (40), (41) and (42) we obtain, for DG

$$\lim_{N\to+\infty} \frac{E_r(\theta)}{\frac{N^2\theta}{2}(\ln N+\gamma)} = 1+e^{-\beta|\theta-c|}.$$

and for PGG

$$\lim_{N\to+\infty} \frac{E_r(\theta)}{\frac{N^2\theta}{2}(\ln N+\gamma)} = 1+e^{-\beta|\theta-c|}e^{\beta c\frac{r}{n}}.$$

$\square$

**Proposition 1.12.** *We have*

*(i) (weak selection limit)*

$$\lim_{\beta\to 0} E_r(\theta) = N^2\left(\sum_{j=1}^{N-1}\frac{1}{j}\right)\theta. \tag{43}$$

*(ii) (large selection limit)*

$$\lim_{\beta\to+\infty} E_r(\theta) = \begin{cases} \frac{N^2}{2}\theta\left(\frac{1}{N-1}+H_N\right) & \text{for } \theta < -\delta, \\ N^2\theta H_N & \text{for } \theta = -\delta, \\ \frac{N^2}{2}\theta\left(1+H_N\right) & \text{for } \theta > -\delta. \end{cases} \tag{44}$$

26

*Proof.* The total expected cost for institutional reward can be rewritten as

$$E_r(\theta) = \frac{N^2\theta(1+e^x)}{2(1+e^x+\ldots+e^{(N-1)x})}\left[\left(1+e^x+\ldots+e^{(N-2)x}\right)\sum_{j=1}^{N-1}\frac{1}{j} + \sum_{j=1}^{N-1}\frac{e^{(N-1-j)x}}{j}\right]$$

$$= \frac{N^2\theta}{2}\frac{\sum_{j=0}^{N-1}\eta_j e^{jx}}{\sum_{j=0}^{N-1}e^{jx}},$$

where

$$\eta_0 := \frac{1}{N-1} + H_N, \quad \eta_j := \frac{1}{N-j-1} + 2H_N + \frac{1}{N-j}, \quad j=1,\ldots,N-2 \quad \text{and} \quad \eta_{N-1} = 1+H_N.$$

We recall that $x = \beta(\theta+\delta)$. Since $\beta \to 0$ implies $x \to 0$, it follows from the formula of $E_r(\theta)$ that

$$\lim_{\beta\to 0} E_r(\theta) = N^2\Big(\sum_{j=1}^{N-1}\frac{1}{j}\Big)\theta.$$

Now we consider the limit $\beta \to +\infty$. If $\theta = -\delta$, then $x = 0$. Hence we also have

$$\lim_{\beta\to+\infty} E_r(\theta) = N^2\theta\frac{\sum_{j=0}^{N-1}\eta_j}{N} = N^2\theta H_N.$$

If $\theta > -\delta$, then $\lim_{\beta\to+\infty}\left(e^{-j(\theta+\delta)}\right)^\beta = 0$ for all $1 \leq j \leq N-2$. Hence

$$\lim_{\beta\to+\infty} E_r(\theta) = \frac{N^2\theta}{2}\lim_{\beta\to+\infty}\frac{\sum_{j=0}^{N-1}\eta_j\left(e^{j(\theta+\delta)}\right)^\beta}{\sum_{j=0}^{N-1}\left(e^{j(\theta+\delta)}\right)^\beta}$$

$$= \frac{N^2\theta}{2}\lim_{\beta\to+\infty}\frac{\sum_{j=0}^{N-1}\eta_j\left(e^{(j-N+1)(\theta+\delta)}\right)^\beta}{\sum_{j=0}^{N-1}\left(e^{(j-N+1)(\theta+\delta)}\right)^\beta}$$

$$= \frac{N^2\theta}{2}\frac{\eta_{N-1}}{1} = \frac{N^2\theta}{2}(1+H_N).$$

If $\theta < -\delta$, then $\lim_{\beta\to+\infty}\left(e^{j(\theta+\delta)}\right)^\beta = 0$ for all $j = 1,\ldots,N-1$. Hence

$$\lim_{\beta\to+\infty} E_r(\theta) = \frac{N^2\theta}{2}\lim_{\beta\to+\infty}\frac{\sum_{j=0}^{N-1}\eta_j\left(e^{j(\theta+\delta)}\right)^\beta}{\sum_{j=0}^{N-1}\left(e^{j(\theta+\delta)}\right)^\beta}$$

$$= \frac{N^2\theta}{2}\frac{\eta_0}{1} = \frac{N^2}{2}\theta\Big(\frac{1}{N-1}+H_N\Big).$$

□

## 1.8 Phase transition

In this section, we prove Theorem 3.2 for reward incentive.

*Proof of Theorem 3.2 for reward incentive.* First we show the existence of phase transition (i.e., the existence of the threshold value $\beta^*$) and study the optimization problem

$$\min_{\theta \geq \theta_0} E_r(\theta), \tag{45}$$

when the selection is less than or equal the threshold value, $\beta \leq \beta^*$. According to Lemma 1.3 we have

$$E'_r(\theta) = \frac{1}{g(x)^2}\left[f(x)g(x) - (x+\beta a)(f(x)g'(x) - f'(x)g(x))\right] = \frac{1}{g(x)^2}\left[Q(u) - (\log u + \beta a)uP(u)\right],$$

where $u = e^x$, $f(x), g(x), P(u)$ and $Q(u)$ are given in Lemma 1.3. By Proposition 1.6, $P$, which is a polynomial with positive leading coefficient, has a unique positive root $u_0 > 1$. Therefore for $0 < u \leq u_0$, $uP(u) < 0$ and hence $E'_r(\theta) \geq 0$ (since $Q(u) > 0$ and $x + \beta a > 0$). Consider now the case $u > u_0$ (thus $uP(u) > 0$). We write $E'_r(\theta)$ as

$$E'_r(\theta) = \frac{1}{g(x)^2} uP(u)\left[F(u) - \beta a\right],$$

where $F(u) = \frac{Q(u)}{uP(u)} - \log u$, which is also given in Lemma 1.3. By Proposition 1.8, $F(u)$ has a global minimizer on $(u_0, +\infty)$. Let $F^*$ be the global minimum. Let $\beta^* := \frac{F^*}{a}$. Then for $\beta \leq \beta^*$, $\beta a \leq \beta^* a = F^* \leq F(u)$ for all $u > u_0$. Thus $F(u) - \beta a \leq 0$ for all $u > u_0$, which implies that $E'(\theta) \geq 0$ for all $u > u_0$ as well. Hence, for $\beta \leq \beta^*$, $E'(\theta) \geq 0$ for all $\theta$ and so, $\min_{\theta \geq \theta_0} E_r(\theta) = E_r(\theta_0)$ and the minimum is attained at $\theta = \theta_0$.

Now we study the optimization problem (45) when the selection is greater than the threshold value, $\beta > \beta^*$. According to Proposition 1.8, $F'$ changes signs at least once when $u > u_0$. Thus $F(u)$ is non-monotonic and the equation $F(u) = \beta a$ has at least two roots which are larger than $u_0$ for all $\beta > \beta^*$. Hence $F(u) - \beta a$, thus $E'_r(\theta)$, changes sign at least twice for any $\beta > \beta^*$. For $N \leq N_0 = 100$, the equation $F(u) = \beta a$ has exactly two solutions $u_0 < u_1 < u_2$. Hence

$$F(u) - \beta a \begin{cases} > 0 & u > u_1, \\ < 0 & u_1 < u < u_2, \\ > 0 & u > u_2. \end{cases}$$

The statement of the main theorem then follows by setting $\theta_i := \frac{\log u_i}{\beta} + a$, $i = 1, 2$. □

## 2 Detailed computations: punishment incentive

In this section, we provide detailed analysis for the case of punishment. We will only sketch on the parts that are similar to the reward case and focus on the parts that are distinguish from the reward case.

## 2.1 The cost function and its derivative

Similarly as in the reward case, using Lemma 1.2 we obtain the expected total cost of interference for institutional punishment

$$E_p(\theta) = \frac{\theta}{2} \sum_{i=1}^{N-1} (n_{1,i} + n_{N-1,i})(N-i)$$

$$= \frac{1}{2} \sum_{j=1}^{N-1} \frac{N^2}{(N-j)j} \Big((W^{-1})_{1j} + (W^{-1})_{N-1,j}\Big)(N-j)\theta$$

$$= \frac{N^2\theta(1+e^x)}{2(e^{Nx}-1)} \sum_{j=1}^{N-1} \frac{e^{(N-1)x} - e^{(j-1)x} + e^{jx} - 1}{j}$$

$$= \frac{N^2\theta(1+e^x)}{2(e^{Nx}-1)} \left[ \big(e^{(N-1)x} - 1\big) \sum_{j=1}^{N-1} \frac{1}{j} + (e^x - 1) \sum_{j=1}^{N-1} \frac{e^{(j-1)x}}{j} \right]$$

$$= \frac{N^2\theta(1+e^x)}{2(1+e^x+\ldots+e^{(N-1)x})} \left[ \big(1 + e^x + \ldots + e^{(N-2)x}\big) \sum_{j=1}^{N-1} \frac{1}{j} + \sum_{j=1}^{N-1} \frac{e^{(j-1)x}}{j} \right]$$

$$:= \frac{N^2\theta}{2} \frac{\hat{f}(x)}{g(x)}, \tag{46}$$

where $g(x) = 1 + e^x + \ldots + e^{(N-1)x}$ is the same as in the reward case, and

$$\hat{f}(x) = (1+e^x)\left[\big(1 + e^x + \ldots + e^{(N-2)x}\big)H_N + \sum_{j=1}^{N-1} \frac{e^{(j-1)x}}{j}\right].$$

The derivative of $E_p$ with respect to $\theta$ is

$$E'_p(\theta) = \frac{N^2}{2} \frac{1}{g(x)^2} \Big[\hat{f}(x)g(x) - (x+\beta a)(\hat{f}(x)g'(x) - \hat{f}'(x)g(x))\Big]. \tag{47}$$

Using again the notation $u = e^x$, by similar computations as in the reward case, we have

$$\hat{f}(x)g(x) - \hat{f}'(x)g(x) = u\hat{P}(u), \quad \hat{f}(x)g(x) = \hat{Q}(u),$$

where

$$\hat{P}(u) := (1+u)\left[\Big(\sum_{j=0}^{N-2}(H_N + \frac{1}{1+j})u^j\Big)\Big(\sum_{j=1}^{N-1} ju^{j-1}\Big) - \Big(\sum_{j=1}^{N-2}\Big(H_N + \frac{1}{1+j}\Big)ju^{j-1}\Big)\Big(\sum_{j=0}^{N-1} u^j\Big)\right]$$

$$- \Big(\sum_{j=0}^{N-2}\Big(H_N + \frac{1}{1+j}\Big)u^j\Big)\Big(\sum_{j=0}^{N-1} u^j\Big) \tag{48}$$

and

$$\hat{Q}(u) := (1+u)\Big(\sum_{j=0}^{N-2}(H_N + \frac{1}{j+1})u^j\Big)\Big(\sum_{j=0}^{N-1} u^j\Big). \tag{49}$$

If $\hat{f}(x)g(x) - \hat{f}'(x)g(x) \leq 0$ then $E'_p(\theta) > 0$. We consider the case $\hat{f}(x)g(x) - \hat{f}'(x)g(x) > 0$, and write $E'_p(\theta)$ as

$$E'_p(\theta) = \frac{N^2}{2} \frac{1}{g(x)^2} (\hat{f}(x)g(x) - \hat{f}'(x)g(x))\Big(\hat{F}(u) - \beta a\Big), \tag{50}$$

where

$$\hat{F}(u) := \frac{\hat{Q}(u)}{u\hat{P}(u)} - \log(u). \tag{51}$$

## 2.2 Concrete examples of small populations

## 2.3 The difference between the reward and punishment costs

In this section, we calculate the difference between the reward and punishment costs.

**Proposition 2.1.**
$$(E_r - E_p)(\theta) = \begin{cases} < 0, & \text{for } \theta < -\delta, \\ = 0, & \text{for } \theta = -\delta, \\ > 0, & \text{for } \theta > -\delta. \end{cases} \tag{52}$$

*Proof.* We have

$$(E_r - E_p)(\theta) = \frac{N^2\theta(1+e^x)}{2(1+e^x+\ldots+e^{(N-1)x})} \sum_{j=1}^{N-1} \frac{e^{(N-1-j)x} - e^{(j-1)x}}{j}$$

$$= \frac{N^2\theta(1+e^x)}{2(1+e^x+\ldots+e^{(N-1)x})} \sum_{j=1}^{N-1} e^{(j-1)x} \left[\frac{1}{N-j} - \frac{1}{j}\right].$$

Hence the sign of $(E_r - E_p)(\theta)$ depends on the sign of the function

$$r(u) := \sum_{j=1}^{N-1} u^{(j-1)} \left[\frac{1}{N-j} - \frac{1}{j}\right]$$

for $u := e^x \in (0, \infty)$. $r(u)$ is a polynomial of degree $N-2$ with a positive leading coefficient. We have

$$r(0) = \frac{1}{N-1} - 1 < 0, \quad \lim_{u \to +\infty} r(u) = +\infty.$$

Moreover, the sequence of coefficients of $r$ has one change of signs, which happens at $j = \lfloor N/2 \rfloor$. Thus, by Decartes's rule of signs, $r$ has a unique positive solution, which is $u = 1$. Hence $r(u) < 0$ for $u \in (0, 1)$ while $r(u) > 0$ for $u \in (1, +\infty)$. As a consequence, $E_r(\theta) < E_p(\theta)$ for $u < 1$, that is for $\theta < -\delta$ and $E_r(\theta) > E_p(\theta)$ for $\theta > -\delta$. □

## 2.4 The polynomial $\hat{P}$

The following proposition establishes interesting relationships between the polynomials $P$ and $\hat{P}$ as well as between $Q$ and $\hat{Q}$. These relationships enable us to obtain properties of the polynomials $\hat{P}$ and $\hat{Q}$ from $P$ and $Q$ without the need to calculate their coefficients, which significantly reduces the complexity of the analysis.

**Proposition 2.2.** *We have*

$$P(u) = -u^{2N-4}\hat{P}\left(\frac{1}{u}\right) \quad \text{and} \quad Q(u) = u^{2N-2}\hat{Q}\left(\frac{1}{u}\right) \tag{53}$$

*or equivalently*

$$\hat{P}(u) = -u^{2N-4}P\left(\frac{1}{u}\right) \quad \text{and} \quad \hat{Q}(u) = u^{2N-2}Q\left(\frac{1}{u}\right). \tag{54}$$

*As a consequence, since $P$ has a unique positive root according to Proposition 1.6, $\hat{P}$ also has a unique positive root.*

*Proof.* By change of index $\hat{j} := N - 2 - j$, we have

(i)

$$\sum_{j=0}^{N-2}\left(H_N+\frac{1}{N-1-j}\right)u^j = \sum_{\hat{j}=0}^{N-2}\left(H_N+\frac{1}{\hat{j}+1}\right)u^{N-2-\hat{j}}$$

$$= u^{N-2}\sum_{\hat{j}=0}^{N-2}\left(H_N+\frac{1}{\hat{j}+1}\right)\left(\frac{1}{u}\right)^{\hat{j}}$$

(ii)

$$\sum_{j=0}^{N-2}\left(H_N+\frac{1}{N-1-j}\right)ju^{j-1} = \sum_{\hat{j}=0}^{N-2}\left(H_N+\frac{1}{\hat{j}+1}\right)(N-2-\hat{j})u^{N-3-\hat{j}}$$

$$= u^{N-4}\sum_{\hat{j}=0}^{N-2}\left(H_N+\frac{1}{\hat{j}+1}\right)(N-2-\hat{j})\left(\frac{1}{u}\right)^{\hat{j}-1}.$$

By change of index $\hat{j} := N-1-j$ we have

(a)

$$\sum_{j=0}^{N-1} ju^{j-1} = \sum_{\hat{j}=0}^{N-1}(N-1-\hat{j})u^{N-2-\hat{j}} = u^{N-3}\sum_{\hat{j}=0}^{N-1}(N-1-\hat{j})\left(\frac{1}{u}\right)^{\hat{j}-1}$$

(b)

$$\sum_{j=0}^{N-1} u^j = \sum_{\hat{j}=0}^{N-1} u^{N-1-\hat{j}} = u^{N-1}\sum_{\hat{j}=0}^{N-1}\left(\frac{1}{u}\right)^{\hat{j}}$$

Therefore, by denoting

$$v := \frac{1}{u}, \quad \hat{a}_j := H_N + \frac{1}{\hat{j}+1}$$

$$P(u) = (1+u)\left[\left(\sum_{j=0}^{N-2}(H_N + \frac{1}{N-1-j})u^j\right)\left(\sum_{j=1}^{N-1} ju^{j-1}\right) - \left(\sum_{j=1}^{N-2}\left(H_N + \frac{1}{N-1-j}\right)ju^{j-1}\right)\left(\sum_{j=0}^{N-1} u^j\right)\right]$$
$$- \left(\sum_{j=0}^{N-2}\left(H_N + \frac{1}{N-1-j}\right)u^j\right)\left(\sum_{j=0}^{N-1} u^j\right)$$
$$= u^{2N-4}\left(1+\frac{1}{u}\right)\left[\left(\sum_{\hat{j}=0}^{N-2}\left(H_N + \frac{1}{\hat{j}+1}\right)\left(\frac{1}{u}\right)^{\hat{j}}\right)\left(\sum_{\hat{j}=0}^{N-1}(N-1-\hat{j})\left(\frac{1}{u}\right)^{\hat{j}-1}\right)\right.$$
$$\left. - \left(\sum_{\hat{j}=0}^{N-2}\left(H_N + \frac{1}{\hat{j}+1}\right)(N-2-\hat{j})\left(\frac{1}{u}\right)^{\hat{j}-1}\right)\left(\sum_{\hat{j}=0}^{N-1}\left(\frac{1}{u}\right)^{\hat{j}}\right)\right]$$
$$- u^{2N-3}\left(\sum_{\hat{j}=0}^{N-2}\left(H_N + \frac{1}{\hat{j}+1}\right)\left(\frac{1}{u}\right)^{\hat{j}}\right)\left(\sum_{\hat{j}=0}^{N-1}\left(\frac{1}{u}\right)^{\hat{j}}\right)$$
$$= u^{2N-4}\left\{(1+v)\left(\left(\sum_{j=0}^{N-2}\hat{a}_j v^j\right)\left[(N-1)\sum_{j=0}^{N-1} v^{j-1} - \sum_{j=0}^{N-1} jv^{j-1}\right]\right.\right.$$
$$\left.\left. - \left[(N-2)\sum_{j=0}^{N-2}\hat{a}_j v^{j-1} - \sum_{j=0}^{N-2} a_j j v^{j-1}\right]\left(\sum_{j=0}^{N-1} v^j\right)\right) - \left(\sum_{j=0}^{N-2}\hat{a}_j v^j\right)\left(\sum_{j=0}^{N-1} v^{j-1}\right)\right\}$$
$$= u^{2N-4}\left\{(1+v)\left[\left(\sum_{j=0}^{N-2}\hat{a}_j v^j\right)\left(\sum_{j=0}^{N-1} v^{j-1}\right) - \left(\sum_{j=0}^{N-2}\hat{a}_j v^j\right)\left(\sum_{j=0}^{N-2} jv^{j-1}\right) + \left(\sum_{j=0}^{N-2}\hat{a}_j j v^{j-1}\right)\left(\sum_{j=0}^{N-1} v^j\right)\right]\right.$$
$$\left. - \left(\sum_{j=0}^{N-2}\hat{a}_j v^j\right)\left(\sum_{j=0}^{N-1} v^{j-1}\right)\right\}$$
$$= u^{2N-4}\left\{(1+v)\left[-\left(\sum_{j=0}^{N-2}\hat{a}_j v^j\right)\left(\sum_{j=0}^{N-2} jv^{j-1}\right) + \left(\sum_{j=0}^{N-2}\hat{a}_j j v^{j-1}\right)\left(\sum_{j=0}^{N-1} v^j\right)\right]\right.$$
$$\left. + \left(\sum_{j=0}^{N-2}\hat{a}_j v^j\right)\left(\sum_{j=0}^{N-1} v^j\right)\right\}$$
$$= -u^{2N-4}\hat{P}(v)$$
$$= -u^{2N-4}\hat{P}\left(\frac{1}{u}\right).$$

Similarly for $Q(u)$, we obtain

$$Q(u) = (1+u)\left(\sum_{j=0}^{N-2}(H_N + \frac{1}{N-1-j})u^j\right)\left(\sum_{j=0}^{N-1} u^j\right)$$
$$= u^{2N-2}\left(1+\frac{1}{u}\right)\left(\sum_{\hat{j}=0}^{N-2}\left(H_N + \frac{1}{\hat{j}+1}\right)\left(\frac{1}{u}\right)^{\hat{j}}\right)\left(\sum_{\hat{j}=0}^{N-1}\left(\frac{1}{u}\right)^{\hat{j}}\right)$$
$$= u^{2N-2}\hat{Q}\left(\frac{1}{u}\right).$$

□

The following lemma, which is the counter-part of Lemma 1.7, establishes a relation between $\hat{P}$ and $\hat{Q}$. It is interesting to notice the appearance of the factor $\frac{3}{4(N-1)}$ on the right-hand side of (55) below.

**Lemma 2.3.** *We have*
$$\hat{Q}(u) \geq \frac{3}{4(N-1)}(1+u)u\hat{P}(u). \qquad (55)$$

*Proof.* It follows from Proposition 2.2 that $\hat{p}_k = -p_{2N-4-k}$. Hence from Proposition 1.6 we obtain

$$\hat{p}_k = \begin{cases} -(H_N + \frac{1}{2}), & k = 0, \\ -(4k+7) + (k+3)H_N - (2k+4)\sum_{j=k+3}^{N-1}\frac{1}{j} + \frac{1}{k+2} & 1 \leq k \leq N-4, \\ -4N+5+\frac{1}{N-1}+NH_N, & k = N-3, \\ 2NH_N - 4(N-1), & k = N-2, \\ -8N+4(k+3)+(2k+4)\sum_{j=k+3-N}^{N-1}\frac{1}{j} + \frac{2N-3-k}{k+2-N} + (2N-k-3)H_N, & N-1 \leq k \leq 2N-5, \\ H_N + \frac{1}{2}, & k = 2N-4. \end{cases}$$

The statement (55) is equivalent to

$$\sum_{j=0}^{N-2}\left(H_N + \frac{1}{1+j}\right)u^j\left(\sum_{j=0}^{N-1}u^j\right) = \sum_{k=0}^{2N-3}\left(\sum_{j=\max(0,k-N+1)}^{\min(N-2,k)}\hat{a}_j\right)u^k \geq \frac{3}{4(N-1)}u\hat{P}(u), \qquad (56)$$

which we now prove. Let

$$\bar{h}_k := \sum_{j=\max(0,k-N+1)}^{\min(N-2,k)}\hat{a}_j = \begin{cases} \sum_{j=0}^{k}\hat{a}_j & 0 \leq k \leq N-2, \\ \sum_{j=k-N+1}^{N-2}\hat{a}_j & N-1 \leq k \leq 2N-3 \end{cases}$$

$$= \begin{cases} (k+1)H_N + \sum_{j=0}^{k}\frac{1}{j+1} & 0 \leq k \leq N-2, \\ (2N-k-2)H_N + \sum_{j=k-N+1}^{N-2}\frac{1}{j+1} & N-1 \leq k \leq 2N-3 \end{cases}$$

To prove (56) we will prove a slightly stronger statement that $\bar{h}_{k+1} \geq \frac{1}{2}\hat{p}_k$ for $k = 0, \ldots, N-2$ and $k = 2N-4$, and that $\bar{h}_{k+1} \geq \frac{3}{4(N-1)}\hat{p}_k$ for $k = N-1, \ldots, 2N-5$.

For $k = 0$: $\bar{h}_1 = 2H_N + 1 + 3/2 > -\frac{1}{2}(H_N + 1/2)$.

For $1 \leq k \leq N-4$, we have

$$\bar{h}_{k+1} - \frac{1}{2}\hat{p}_k = (k+2)H_N + \sum_{j=0}^{k+1}\frac{1}{j+1} + \frac{1}{2}\left((4k+7) - (k+3)H_N + (2k+4)\sum_{j=k+3}^{N-1}\frac{1}{j} - \frac{1}{k+2}\right)$$

$$= (k+1)(H_N - H_{k+3}) + \frac{(4k+7)}{2} + \frac{k+3}{2}H_N - \frac{1}{2(k+2)} > 0.$$

For $k = N-3$:

$$\bar{h}_{k+1} - \frac{1}{2}\hat{p}_k = (N-1)H_N + \sum_{j=0}^{N-2}\frac{1}{1+j} - \frac{1}{2}(-4N+5+\frac{1}{N-1}+NH_N)$$

$$= \frac{1}{2}\left((4N-5) + NH_N - \frac{1}{N-1}\right) > 0.$$

For $k = N-2$

$$\bar{h}_{k+1} - \frac{1}{2}\hat{p}_k = (N-1)H_N + \sum_{j=0}^{N-2}\frac{1}{1+j} - \frac{1}{2}(2NH_N - 4(N-1))$$

$$= 2(N-1) > 0.$$

33

For $N - 1 \leq k \leq 2N - 5$. Let $\alpha > 1$ be chosen later. We have

$$\bar{h}_{k+1} - \frac{1}{\alpha}\hat{p}_k = (2N - k - 3)H_N + \sum_{j=k+3-N}^{N-1} \frac{1}{j}$$

$$- \frac{1}{\alpha}\Big(-8N + 4(k+3) + (2k+4)\sum_{j=k+3-N}^{N-1}\frac{1}{j} + \frac{2N-3-k}{k+2-N} + (2N-k-3)H_N\Big)$$

$$= (2N - k - 3)H_N\Big(1 - \frac{1}{\alpha}\Big) + \Big(1 - \frac{2(k+2)}{\alpha}\Big)\sum_{j=k+3-N}^{N-1}\frac{1}{j} - \frac{(2N-3-k)}{\alpha(k+2-N)}$$

$$\geq (2N - k - 3)\Big(1 - \frac{1}{\alpha}\Big)H_N - \Big(\frac{2(k+2)}{\alpha} - 1\Big)\sum_{j=k+2-N}^{N-1}\frac{1}{j}.$$

We choose $\alpha > 1$ such that

$$(2N - k - 3)\Big(1 - \frac{1}{\alpha}\Big) \geq \Big(\frac{2(k+2)}{\alpha} - 1\Big) \quad \text{for all} \quad k = N - 1, \ldots, 2N - 5.$$

The above condition is equivalent to the following condition

$$2N(\alpha - 1) \geq (k+2)\alpha + (k+1) \quad \text{for all} \quad k = N - 1, \ldots, 2N - 5.$$

Therefore we only need to choose $\alpha > 1$ such that the above inequality holds true for $k = 2N - 5$, that is

$$2N(\alpha - 1) \geq (2N - 3)\alpha + (2N - 4).$$

Hence $\alpha \geq \frac{4(N-1)}{3}$. Therefore, we obtain that for $N - 1 \leq k \leq 2N - 5$

$$\bar{h}_{k+1} \geq \frac{3}{4(N-1)}\hat{p}_k.$$

For $k = 2N - 4$, we have

$$\bar{h}_{k+1} - \frac{1}{2}\hat{p}_k = H_N + \frac{1}{N-1} - \frac{1}{2}(H_N + \frac{1}{2}) = \frac{1}{2}H_N + \frac{1}{N-1} - \frac{1}{4} > 0.$$

$\square$

Let $\hat{u}^0$ be the unique positive root of $\hat{P}$. Note that $\hat{u}^0 = \frac{1}{u^0}$, where $u^0$ is the unique positive root of $P$. We consider the minimisation problem

$$\min\{\hat{F}(u) : \quad u > \hat{u}^0\}, \tag{57}$$

where $\hat{F}$ is defined in (51). The following proposition is the counter-part of Proposition 1.8.

**Proposition 2.4.** *The minimisation problem (57) has a global minimizer. For $N \leq N_0 = 100$ sufficiently small, it has a unique minimizer $\hat{u}^*$ and $\hat{F}$ is decreasing when $\hat{u}^0 < u < \hat{u}^*$ and is increasing when $u > \hat{u}^*$.*

*Proof.* The proof of this proposition is similar to that of Proposition 1.8 in the reward case. According to Lemma 2.3 we have

$$\hat{F}(u) \geq \frac{3}{4(N-1)}(1+u) - \log(u).$$

It implies that

$$\lim_{u \to +\infty} \hat{F}(u) = +\infty = \lim_{u \to \hat{u}_0^+} \hat{F}(u).$$

From this together with the fact that $\hat{F}$ is smooth on the interval $(\hat{u}^0, +\infty)$, we conclude that $\hat{F}$ has a minimizer on this domain. The rest of the proof is the same as that of Proposition 1.8 where $P$ and $Q$ are replaced by $\hat{P}$ and $\hat{Q}$, respectively and the unique positive root $u^0$ of $P$ is replaced by that of $\hat{P}$ (i.e., $\hat{u}^0$).

$\square$

## 2.5 Finite population estimates

**Lemma 2.5.** *It holds that*

$$\frac{N^2\theta}{2}\Big(H_N + \frac{1}{N-1}\Big) \leq E_p(\theta) \leq N(N-1)\theta\Big(H_N + 1\Big). \tag{58}$$

*Proof.* We recall that

$$E_p(\theta) = \frac{N^2\theta(1+e^x)}{2(1+e^x+\ldots+e^{(N-1)x})}\bigg[(1+e^x+\ldots+e^{(N-2)x})\sum_{j=1}^{N-1}\frac{1}{j} + \sum_{j=1}^{N-1}\frac{e^{(j-1)x}}{j}\bigg]. \tag{59}$$

Since we also have the following estimates

$$\frac{1}{N-1}(1+e^x+\ldots+e^{(N-2)x}) \leq \sum_{j=1}^{N-1}\frac{e^{(j-1)x}}{j} \leq (1+e^x+\ldots+e^{(N-2)x}),$$

we can proceed exactly as in the proof of Proposition 1.11 in the reward case to obtain the required estimates. □

## 2.6 Asymptotic limits

**Proposition 2.6.** *We have the following limits:*

(i) *(Infinite-population limit)*

$$\lim_{N\to+\infty}\frac{E_p(\theta)}{\frac{N^2\theta}{2}(\ln N + \gamma)} = \begin{cases} 1 + e^{-\beta|\theta-c|} & \text{for } DG, \\ 1 + e^{-\beta|\theta-c|}e^{\beta c\frac{r}{n}} & \text{for } PGG. \end{cases}$$

(ii) *(weak selection limit)*

$$\lim_{\beta\to 0} E_p[(\theta) = N^2\theta H_N.$$

(iii) *(large selection limit)*

$$\lim_{\beta\to+\infty} E_p(\theta) = \begin{cases} \frac{N^2\theta}{2}\Big(1 + H_N\Big) & \text{for } \theta < -\delta, \\ N^2\theta H_N & \text{for } \theta = -\delta, \\ \frac{N^2\theta}{2}\Big(H_N + \frac{1}{N-1}\Big) & \text{for } \theta > -\delta. \end{cases}$$

*Proof.*

(i) We also have

$$\frac{N^2\theta}{2}\frac{(1+e^x)(1+\ldots+e^{(N-2)x})}{1+\ldots+e^{(N-1)x}}\Big(\sum_{j=1}^{N-1}\frac{1}{j} + \frac{1}{N-1}\Big) \leq E_p(\theta)$$

$$\leq \frac{N^2\theta}{2}\frac{(1+e^x)(1+\ldots+e^{(N-2)x})}{1+\ldots+e^{(N-1)x}}\Big(\sum_{j=1}^{N-1}\frac{1}{j} + 1\Big)$$

Hence by proceeding exactly the same as in the reward case, we obtain the desired limit.

(ii) Since $\beta \to 0^+$ implies $x \to 0$, we have

$$\lim_{\beta\to 0^+} E_p(\theta) = N^2\theta\Big(\sum_{j=1}^{N-1}\frac{1}{j}\Big).$$

(iii) We have

$$E_p(\theta) = \frac{N^2\theta(1+e^x)}{2(1+e^x+\ldots+e^{(N-1)x})}\left[(1+e^x+\ldots+e^{(N-2)x})\sum_{j=1}^{N-1}\frac{1}{j}+\sum_{j=1}^{N-1}\frac{e^{(j-1)x}}{j}\right].$$

$$= \frac{N^2\theta}{2}\frac{\sum_{j=0}^{N-1}\hat{\eta}_j e^{jx}}{\sum_{j=0}^{N-1}e^{jx}},$$

where

$$\hat{\eta}_0 = H_N + 1, \quad \hat{\eta}_j = 2H_N + \frac{1}{j} + \frac{1}{1+j} \quad \text{for } j=1,\ldots,N-2 \quad \text{and} \quad \hat{\eta}_{N-1} = H_N + \frac{1}{N-1}.$$

Proceeding similarly as in Proposition 1.12 we get

$$\lim_{\beta\to+\infty} E_p(\theta) = \begin{cases} \frac{N^2\theta}{2}\frac{\hat{\eta}_0}{1} = \frac{N^2\theta}{2}\left(1+H_N\right) & \text{for } \theta < -\delta, \\ N^2\theta H_N & \text{for } \theta = -\delta, \\ \frac{N^2\theta}{2}\frac{\hat{\eta}_{N-1}}{1} = \frac{N^2\theta}{2}\left(H_N+\frac{1}{N-1}\right) & \text{for } \theta > -\delta. \end{cases}$$

□

## 2.7 Phase transition

Having obtained Proposition 2.4, the proof of Theorem 3.2 for the case of punishment can be proceeded in the exactly as in the reward case by replacing $F^*$ by $\hat{F}^*$, which is the minimum value of $\hat{F}$ in the interval $(\hat{u}^0, +\infty)$ and $\beta^*$ by $\beta_p^*$, where

$$\beta_p^* := \frac{\hat{F}^*}{a}. \tag{60}$$

Therefore we omit the details.


%
%
%
\end{document}